\documentclass[10pt,a4paper]{article} 
\usepackage{slashbox} 
\usepackage[latin1]{inputenc} 
\usepackage[dvips]{graphicx} 
\usepackage{epsfig} 
\usepackage{float} 
\usepackage{graphicx} 
\usepackage[T1]{fontenc} 
\usepackage{indentfirst} 
\usepackage{enumerate} 
\usepackage{amssymb} 
\usepackage{amsmath} 
 
\newtheorem{remarque}{Remark}[section] 
 
\def\AC{{\cal A}} 
\def\TC{{\cal T}} 
\def\NC{{\cal N}} 
\def\HC{{\cal H}} 
\def\EC{{\cal E}} 
\def\UC{{\cal U}} 
\def\RB{{\Bbb R}} 
\def\RC{{\cal R}}

\newcommand{\be}{\begin{eqnarray}} 
\newcommand{\ee}[1]{\label{eq:#1}\end{eqnarray}} 
 
\newcommand{\ese}{\end{eqnarray*}} 
\newcommand{\bse}{\begin{eqnarray*}}

\newtheorem{col}{Corollary} 
\newcommand{\bcol}{\begin{col}} 
\newcommand{\ecol}[1]{\label{the:#1}\end{col}\par} 
 
\newtheorem{defi}{Definition} 
\newcommand{\bdf}{\begin{defi}} 
\newcommand{\edf}[1]{\label{df:#1}\end{defi}\par} 
 
\newtheorem{lem}{Lemma} 
\newcommand{\blem}{\begin{lem}} 
\newcommand{\elem}[1]{\label{le:#1}\end{lem}\par} 
 
\newtheorem{pro}{Proposition} 
\newcommand{\bpro}{\begin{pro}} 
\newcommand{\epro}[1]{\label{pr:#1}\end{pro}\par} 
 
\newtheorem{statement}{Problem}

\newcommand{\card}{\mbox{card}}

\newcommand{\bstatement}{\begin{statement}} 
\newcommand{\estatement}[1]{\label{stat:#1}\end{statement}\par} 
 
\newcounter{assc} 
\newcommand{\bass}[1]{\refstepcounter{assc}\label{ass:#1} \begin{maliste}{\bf \arabic{assc}.}} 
\newcommand{\eass}{\end{maliste}}

\title{ 
{\bf A Value-At-Risk Approach for Robust Management of Electricity Power Generation}} 
 
\author{Vincent GUIGUES\footnote{SMS/LMC, universit\'e Joseph Fourier,
BP 53, F38041 Grenoble Cedex 9. Email: vincent.guigues@imag.fr.} \and
Papa Momar NDIAYE\footnote{Raise Partner SAS. 3 chemin de Ronde, F38000 Grenoble. Email: papa-momar.ndiaye@raisepartner.com.}} 
 
\begin{document} 
 
\date{} 
\maketitle 
 
\section*{Abstract} 
In this paper, we apply Value-At-Risk (VaR) approaches on the problem of yearly electric generation management. In a classical 
approach, the future is modelled as a markov chain and the goal is to minimize the average generation cost over this uncertain future. However, such a strategy could lead to big financial losses if worst case scenarios occur. The two VaR approaches we propose, precisely aim at robustifying the model. On a practical point of view, it amounts to introduce a new set of constraints modelling the uncertainties in the original optimization problem or equivalently to change the dual objective function. The new optimization problems are solved as efficiently as the nominal model. Numerical simulations are presented and discussed for this application. 
 
\section{Introduction} 
 In this paper, we are interested in optimization problems arising in electrical yearly power management. Given electric generation plants (nuclear, thermal and hydroelectric power generator plants, demand side management contracts modelled as a virtual plant called EJP),  the objective is to minimize the production cost over a yearly horizon, to fulfill operating constraints of generation units and the equilibrium between production and demand at each time step. In  practice, the modelling approach is highly depending on the time horizon of the optimization problem : for short time horizons, typically  daily or weekly, the problem is generally assumed to be deterministic (cf. \cite{deter1},\cite{deter2}), whether for longer management horizons,  a special emphasis is done on the stochastic nature of data and events. In particular, on a yearly scale, reservoir inflows, demand, availability of the plants as well as electricity prices cannot be considered to be deterministic : in France, for example, winter costumer's demand has uncertainties that can reach one GW per decreasing temperature degree while the peak loads are around 70 GW ! So a true challenge is to ensure robustness of computed optimal production and marginal value face to various uncertainties like customer demand but also water inflows or plants unavailability. 
 
\medskip 
 
In general, for yearly generation management, utilities are not interested in generation scheduling but rather management strategies or Bellman values to perform Monte Carlo analysis of the futur. Since Stochastic Dynamic Programming  quickly comes to its limits for the optimization of high dimensional state systems, a decomposition approach is usually necessary \cite{spacedecomposition}. In our case, a large scale numerical optimization is first formulated and solved solved using Lagrangian relaxation to provide marginal costs on a scenario tree. Then, those marginal costs are used to compute local feed-backs (see below section 2.1 for details). 
      This decompostion framework that could be compared to the Dual Stochastic Dynamic Programming method \cite{pereiraandpinto} is based on the adaptation of  \cite{carpentier1} to yearly generation management. The main drawback of this scheme relies in the local aspect of the feed-back functions that loses the robustness of the global Bellman function. 
 
In this paper, we show that a Value-At-Risk (VaR) approach for modelling uncertainties allows to enhance the robustness properties of those local feed-back laws with the following practical  benefits:  (i)  significant  reduction of the variance of simulated cost - up to 38\% of reduction for a comparable average cost if the worst case scenario occurs -; (iii)  parsimonious use of water reservoir; (ii)  reduction of very high cost strategies.  Moreover, we will show that  there exits VaR approaches on the dual optimization problem that preserve the space decomposition approach while having a very nice economical and physical interpretation. The first one comes from a primal relaxation on the demand side - how to control the sales turnover face to uncertainty of the demand ?- and the second one from a dual relaxation -how to control  the production  costs face to plant random unavailability ?-.

\medskip 
 
The sequel  is organized as follows. In the next section, the optimization model is described and special focus is given on random events and inputs. Section 3 deals with the robustification issues, the VaR approaches, the connection between VaR and duality as well as implementation issues. Section 4 presents some numerical results comparing the performances of the nominal and robust models and finally some additional modelling details are given in the appendix. 
 
\section{Setting of the physical model} 
The  physical model is a stochastic dynamic system where the random inputs are the costumer's demand, the unavailability of 
the thermal units and the quantity of natural water inflows. With a representation of the random inputs and events as Markov chains, 
we can naturally formulate the optimization problem as a finite horizon discrete time stochastic control problem on a scenario tree representing the behavior of the random inputs and states. 
 
\subsection{Model setting} 
We aim at minimizing the average  production cost along the scenario tree. If we describe at each node  $n$ of the tree the states of the plant production unit  $\ell$ by the variable $x_n^{\ell}$ and the commands applied to this plant by  $u_n^{\ell}$, the stochastic control problem may be formulated as \cite{soprano}: 
\begin{equation}\label{base} 
\left\{ 
\begin{array}{l} 
\displaystyle \min_u \; \sum_{\ell \in \mathcal{L}} \; \sum_{n \in \mathcal{O}} 
\pi_n \mathcal{C}_{n,p}^{\ell} ( x_n^{\ell}(p),u_n^{\ell}(p))\\ 
\displaystyle \forall n \in \mathcal{O}, \; \forall p \in \mathcal{P}_n, \;\sum_{\ell \in \mathcal{L}} P_{n,p}^{\ell} ( x_n^{\ell}(p),u_{n}^{\ell}(p)) = \mathcal{D}_n(p),  \\ 
\forall \ell \in \mathcal{L},  \quad \displaystyle 
(x_{\bullet}^{\ell}(\bullet),u_{\bullet}^{\ell}(\bullet))  \in 
\chi_{_\ell}, 
\end{array} 
\right. 
\end{equation} 
where 
\begin{itemize} 
\item $\mathcal{L}$ is the set of plants and $\mathcal{O}$ is the set of the nodes, 
\item $\pi_n$ is the probability to be at node $n$, 
\item $\mathcal{P}_n$ is the set of time subdivisions associated to node  $n$, 
\item $x_n^{\ell}  (p)$ the state of plant $\ell$ at node $n$ and time subdivision $(p)$, 
\item $u_n^{\ell} (p)$  is the control variable of plant $\ell$ at node $n$ and time subdivision $(p)$, 
\item $P_{n,p}^{\ell} ( x_n^{\ell}(p),u_{n,p}^{\ell}(p))$ is the  production of plant  $\ell$ in the state  $x_n^{\ell}$ when command 
$u_n^{\ell} (p)$ is applied to this plant at node $n$ and time 
subdivision $(p)$, 
\item $\mathcal{D}_n(p)$ is the costumer demand at node $n$ and time subdivision $(p)$, 
\item $C_{n,p}^\ell(x_{n}^\ell(p),u_{n}^\ell(p))$ is the production cost when command $u_{n}^\ell(p)$ is applied to unit $\ell$ in the state $x_{n}^\ell(p)$, 
\item $\chi_{_\ell}$ is the functional set of constraints on the control and state  variables of plant $\ell$. 
\end{itemize} 
 
In fact, due to the autonomy of the plants, the model may be reformulated as a linear optimization problem with separated domains of constraints and one coupling constraint (production/demand equilibrium). Each domain of constraint is a dynamic system describing  a plant process.  Introducing cost vectors  $c_{i}$ and coupling matrices  $\AC_i$  with ad hoc sizes for  $i \in \{1, ..., 4\}$, we can formally describe  \eqref{base} in the following way: 
\begin{equation}\label{base2} 
\left\{ 
\begin{array}{ll} 
\min \;c^T_1u_t + c^T_2 u_n + c^T_3x_h+c^T_4 x_e\\ 
(u_t,u_n,u_h,u_e) \in {\cal T} \times {\cal N} \times {\cal H}(x_h) \times {\cal E}(x_e),\\ 
{\AC}_1u_t + {\AC}_2u_n + {\AC}_3u_h + {\AC}_4u_e   = d \in \RB^D, 
\end{array} 
\right. 
\end{equation} 
where each control variable  $u_{\bullet}$ belongs to a set 
parameterized\footnote{see Appendix A2 for a detailed plants 
description} by the corresponding state variable $x_{\bullet}$ : 
\begin{itemize} 
\item   $u_t$,  control variable of classical thermal plants, with $\cal T$ as the set of constraints for the thermal plants 
subset; 
\item   $u_n$,  control variable of nuclear thermal plants, with $\cal N$ as the set of constraints for the thermal plants  subset; 
\item  $x_h$, state variable of hydraulic plants  with  $u_h$ as the control variable and dynamics described by  $u_h \in {\cal H}(x_h)$; 
\item $x_e$, state variable of  EJP contract, with  $u_e$  as the control variable and dynamics described by   $u_e \in {\cal E}(x_e)$; 
\item $d \in \RB^D$, the vector of demands with $D$ =  (n° of time subdivision ) $\times$ card ($\mathcal{O}$). 
\end{itemize} 
Therefore we may write \eqref{base2} as: 
\begin{equation}\label{base3} 
(LP)\quad \left\{ 
\begin{array}{ll} 
\displaystyle \min_{u \in \UC} f^{\ell}(u,x_h,x_e)\\ 
\AC u  = d \mbox{ given }\in \RB^D, 
\end{array} 
\right. 
\end{equation} 
when setting 
\begin{eqnarray} 
u &=&(u_t,u_n,u_h,u_e) \in \ \UC = {\cal T} \times \ {\NC} \times {\HC}(x_h) \times {\EC}(x_e),\\ 
 f^{\ell}(u, x_h,x_e) &=& \;c^T_1u_t + c^T_2 u_n + c^T_3x_h+c^T_4x_e,\\ 
\AC u &=& {\AC}_1u_t + {\AC}_2u_n + {\AC}_3u_h + {\AC}_4u_e. 
\end{eqnarray} 
Notice that in this model, the demand $d$ is a fixed vector corresponding to the realizations of the demand at the different nodes to the scenario tree.

\subsection{Model Analysis} \label{default} 
 
The efficiency of this model is assessed on a set of independent scenarios representing different evolutions of the demand, the inflows for hydro reservoirs and the outages of the thermal units. For each scenario, a generation schedule as well as its cost are determined (a detailed description of the implementation of the generation schedule is given in the appendix). Such a model has intrinsic limitations essentially linked to the fact that the strategy computed by the algorithm above is optimal only on the optimal trajectory of each reserve. In this sense, the Bellman functions obtained by Stochastic Dynamic Programming on marginal values only give a local optimum. The effect of using such local feedbacks as global strategies is to create a high volatility of scenarios costs. So it is desirable to strengthen this model by including a more reliable model of uncertainty on the scenarios at the earliest stages of the problem setting. The goal is to ensure some regularity of the optimal strategies with respect to the inputs of the optimization problem. In other words we would like to find a robust counterpart with the following properties: 
\begin{itemize} 
\item[(i)] reduce the volatility of the simulated costs over a continuum set of reasonable 
scenarios; 
\item[(ii)] reduce the number of extreme case optimal strategies (parsimonious use of water reservoir that might not be  nearly empty for a long 
period); 
\item[(ii)] reduce the number of very high cost optimal 
strategies. 
\end{itemize} 
 We will see in the next section that such an objective of variance reduction may be easy to formulate in a Value-At-Risk setting. 
 
\section{Robust Counterpart of the Decision Model} 
 
\subsection{The Value-At-Risk approach for stochastic optimization 
problems}\label{var} 
 
Let $(\omega,r,x) \rightarrow f(r(\omega),x)$ be a concave (with respect to $x$) income functional depending on a random  function $\omega \rightarrow r(\omega)$ where  $x \in X \subset \mathbb{R}^n$ is deterministic variable, $X$ being a non empty closed and bounded set. A Value-At-Risk (VaR) approach allows us to choose $x$ leading to the maximal possible income with a given confidence level $0<\varepsilon<1$. Typically, if we have additional constraints on $x$ expressed as $g(x) \geq 0$, one formulates the following optimization problem : 
\begin{equation} \label{vardef} 
\left\{ 
\begin{array}{l} 
\max \gamma \\ 
P(f(r(w),x) \geq \gamma) \geq 1 - \varepsilon,\\ 
x \in X, \quad g(x) \geq 0. 
\end{array} 
\right. 
\end{equation} 
Let $\Phi$ be the cumulative distribution function of the Gaussian 
density. Following \cite{bertsimas} and \cite{smith}, we can find 
a Risk Averse  solution (indeed a prudent one as upper bound of the optimum) in some usual 
cases when solving: 
\begin{equation} \label{vardef} 
(VaR_{\varepsilon}) \left\{ 
\begin{array}{l} 
\displaystyle\max \;\;\;E_{\omega}[f(r(\omega),x)]-\kappa(\varepsilon) \, \sigma_{\omega}[f(r(\omega),x)]\\ 
x \in X, \quad g(x) \geq 0, 
\end{array} 
\right. 
\end{equation} 
where $\kappa(\varepsilon)$ is a risk factor depending on the assumptions on the distribution: 
\begin{equation} \label{kappa} 
\kappa(\varepsilon) = 
\left\{ 
\begin{array}{l} 
{\Phi^{-1}(1-\varepsilon)}>0 \mbox{ if $f(r(\bullet),x)$ is Gaussian},\\ 
\sqrt{\frac{1-\varepsilon}{\varepsilon}} \mbox{ if 
$f(r(\bullet),x) \in L^1_{\RB} \cap  L^2_{\RB}$.} 
\end{array} 
\right. 
\end{equation} 
 For instance is $f(c(\omega),x)=c(\omega)^Tx$ is a linear function of $x$, the problem ($VaR_{\varepsilon}$) simply reduces to 
\begin{equation}  \label{vardef1} 
\left\{ 
\begin{array}{l} 
\displaystyle\max \;\;\;E_{\omega}[c(\omega)^T]x-\kappa(\varepsilon) \sqrt{x^T \Gamma x} \mbox{ with } \Gamma_{ij}=cov(c_i(w),c_j(w)),\\ 
x \in X, \quad g(x) \geq 0, 
\end{array} 
\right. 
\end{equation} 
and if $\Gamma$ is invertible the above problem appears as the Robust Counterpart of the problem 
\begin{equation} 
\left\{ 
\begin{array}{l} 
\displaystyle\max \;\;\;c(\omega)^Tx\\ 
x \in X, \quad g(x) \geq 0, 
\end{array} 
\right. 
\end{equation} 
where the uncertainty set chosen for the random vector $c(\omega)$ 
is the ellipsoid: 
$$ 
\{x \in \mathbb{R}^n,\;\; (x-E_{\omega}[c(\omega)])^T \Gamma^{-1} 
(x-E_{\omega}[c(\omega)]) \leq \kappa^2(\varepsilon)\}. 
$$ 
Now it is clear that a VaR approach is a practical way to 
calibrate  a variance penalty term for a maximization of  a random 
functional. In the case of power generation management, it means 
that we aim to find the best compromise between production cost 
and volatility of strategies at (possibly) the extra cost of some 
sub-optimality on the most favorable scenarios. This approach is 
very easy to set as a regularization approach for linear programs 
because one only needs to specify a risk exposure level 
$\varepsilon$. 
\begin{remarque} {\em Being a covariance matrix,  $\Gamma$ is symmetric and positive semidefinite, and using the 
induced norm, the objective function in \eqref{vardef1}  may be 
written as \begin{math} E_{\omega}[c(\omega)x] - 
\kappa(\varepsilon)\Vert x \Vert_{\Gamma}\end{math}. Therefore a 
VaR regularization appears as an Han Penalization for the problem: 
 \begin{equation} 
\left\{ 
\begin{array}{l} 
\displaystyle\max \;\;\;E_{\omega}[c(\omega,x)]\\ 
x \in X, \quad g(x) \geq 0, \quad \sigma[c(\omega,x)] \leq \alpha, 
\end{array} 
\right. 
\end{equation} 
where the penalization coefficient is chosen on a probabilistic 
basis. We mention that in the convex case, there exists some 
sufficient conditions linking $\kappa(\varepsilon)$ to the dual 
norm of $\Vert \bullet \Vert_{\Gamma}$  that  ensures the 
exactness of this penalization. Remarking that if $\Gamma \succ 0 
$, the dual norm is $\Vert \bullet \Vert_{\Gamma^{-1}}$, we can 
mention that the condition given in \cite{bgls} reduces to 
$\kappa(\varepsilon) \geq \Vert \lambda \Vert _{\Gamma^{-1}}$ 
where $\lambda$ is the lagrange multiplier associated to the 
optimal solution. Therefore : 
\begin{itemize} 
\item [(i)] a sufficient condition for exact penalization is $\varepsilon \leq \frac{1}{1  + \lambda_{min}(\Gamma)}$, where $\lambda_{min}(\Gamma)$ is the smallest eigenvalue of $\Gamma$. 
\item [(ii)] Optimizing the smallest eigenvalue of the covariance matrix $\Gamma$ will allow to reduce the bound 
  on $\varepsilon$ and then to enforce the constraint of risk 
  reduction: $P(c(\omega)x \geq \gamma) \geq 1 -\varepsilon$. So any matrix 
calibration technique that will reduce the 
  condition number of the covariance matrix by increasing the smallest eigenvalue will provide 
  additive degree of freedom if we use a probabilistic constraint to control the income. 
\end{itemize} } 
\end{remarque} 
 
\subsection{Application to the Power Generation Model} 
 
The idea is to take advantage of a decomposition of the dual 
optimization problem to introduce a VaR modelization  on two 
subproblems : first on  the uncertainty on the demand and next on 
the unavailability  of thermoelectric plants.  To begin with, we 
point out that the dual problem of \eqref{base3} is 
\begin{math} \displaystyle{\max_{\lambda \in \mathbb{R}^D}} \theta(\lambda)\end{math} 
where: 
\begin{equation} 
\displaystyle \theta(\lambda) = {\theta}_d(\lambda) + 
{\theta}_T(\lambda)+ {\theta}_N(\lambda)+ \theta_H(\lambda)+ 
\theta_ J(\lambda)= \theta_d(\lambda) + \tilde 
\theta(\lambda),\end{equation} once we have introduced the partial 
dual functions: 
\begin{equation}\label{dualdec} 
\left\{ 
\begin{array}{lll} 
{\theta}_d(\lambda)&=&  \lambda^T d,\\ 
{\theta}_T(\lambda)&=& \displaystyle  \inf_{u_t \in \TC} \; (c_1 - \AC^T_1\lambda)^T \; u_t,\\ 
{\theta}_N(\lambda)&=&  \displaystyle\inf_{u_n \in \NC} \; (c_2 - \AC^T_2\lambda)^T \; u_n,\\ 
{\theta}_H(\lambda)&=& \displaystyle \inf_{u_h \in \HC(x_h)} \; c^T_3x_h -  \lambda^T\AC_3u_h, \\ 
{\theta}_J(\lambda)&=& \displaystyle \inf_{u_e \in \EC(x_e) }  \; 
c^T_4x_e -  \lambda^T \AC_4u_e. 
\end{array}\right. 
\end{equation} 
Now  it is natural to robustify  \eqref{base3} by formulating a 
VaR problem on the subsystem with random information or state : 
${\theta}_d$ for the costumer's demand $d$ and 
$({\theta}_T,{\theta}_N,\theta_H)$ for the states of the thermal, 
nuclear and hydro plants. 
\subsubsection{Primal Relaxation and Dual VaR regularization of the demand} 
The idea is to make  a  primal relaxation  on  the predictions 
made for the demands on the different scenarios (i.e the values of 
the demand $d$ at the different nodes of the scenario tree) that 
are prone to errors. Rather 
  than considering that the demands at each node and each time subdivision of the tree are known exactly, we suppose that $d$  belongs to a given uncertainty set $\mathcal{E}$ which is 
  the ellipsoid given by: 
\begin{equation} 
\mathcal{E} = \mathcal{E}(\bar d, \Gamma, \kappa)=\{\, x \in 
\mathbb{R}^D \;\; / \;\; {(x- \bar d )}^T \Gamma^{-1}(x -\bar d) 
\leq {\kappa}^2(\varepsilon) \,\}. 
\end{equation} 
where  $\bar d = E_{\omega}[d(w)]$, the covariance matrix $\Gamma$ 
is given by  $\Gamma_{ij}=cov(d_i(\omega),d_j\omega))$ and 
$\kappa(\varepsilon)$ depends on the assumptions made on the 
distribution of the demand (see \cite{oustry}). That means that we 
reformulate the problem \eqref{base3} as: 
\begin{equation}\label{base4} 
(LP_{\RC})\quad \left\{ 
\begin{array}{ll} 
\displaystyle \min_{u \in \UC} f^{\ell}(u,x_h,x_e)\\ 
\AC u  = d \in {\EC}. 
\end{array} 
\right. 
\end{equation} 
This is a relaxation of problem \eqref{base3}: the demand vector 
is no longer fixed but can be any vector from the ellipsoid $\mathcal{E}$. 
Solving \eqref{base4} by duality amounts to solve 
\begin{math} \displaystyle 
\max_{\lambda}\theta_{\RC}(\lambda)\end{math} with 
\begin{equation} 
\begin{array}{lll} 
\displaystyle \theta_{\RC}(\lambda)&=& \displaystyle \min_{u \in \UC, \quad d \in {\EC}} f^{\ell}(u,x_h,x_e) + {\lambda}^T(d - {\AC u} )\\&=&  {\theta}_T(\lambda)+ {\theta}_N(\lambda)+ \theta_H(\lambda)+ 
\theta_ J(\lambda) + \displaystyle \min_{d \in {\EC}} {\lambda}^T 
d. 
\end{array} 
\end{equation} 
Now notice that 
\begin{math} 
\displaystyle \min_{d \in {\EC}} {\lambda}^T d = 
\phi_{\EC}(\lambda) 
\end{math} 
where $\phi_{\EC}$ is the support function of the uncertainty set 
$\EC$ and is given by : 
\begin{equation} 
\phi_{\EC}(\lambda) = \lambda^T \bar d - \kappa(\varepsilon) \sqrt 
{\lambda^T \Gamma \lambda}. 
\end{equation} 
Note that the robustification just turns out to replace 
$\theta_d(\lambda) = \lambda^T d $  in the original dual function 
$\theta$ by $\phi_{\EC}(\lambda) = \lambda^T \bar d - \kappa \sqrt 
{\lambda^T \Gamma \lambda}$,  where  $\bar d = 
E_{\omega}[d(\omega)]$. We can notice that the relaxation of the 
demand in the ellispoid $\mathcal{E}$ in problem \eqref{base4} 
amounts to use a VaR approach on the dual problem of problem (LP). 
Indeed, using a VaR orientated technique, as the demand d is 
random, instead of maximizing $\theta(\lambda)$ (which is the dual 
problem of problem (LP)) we could maximize $\tilde \theta(\lambda) 
+ \gamma^*(\lambda)$, where 
\begin{equation}\label{VarST} 
\gamma^*(\lambda)=\left\{ 
\begin{array}{l} 
\max \;\;\gamma\\ 
P(\lambda^T d(\omega) \geq \gamma) \geq 1 - \varepsilon. 
\end{array} 
\right. 
\end{equation} 
From subsection \ref{var} this VaR approach reduces to 
$\displaystyle{\max_{\lambda}} \; \theta_{\mathcal{R}}(\lambda)$, 
which is the dual problem of problem $(LP_{\mathcal{R}})$. In what 
follows, this VaR approach will be denoted by $VaR_{FA}$.\\\\ {\bf 
Economical Interpretation.} The problem \eqref{VarST} may be 
interpreted as the maximization of the minimal Sales Turnover that 
can be ensured with an arbritary degree of confidence. In other 
words, a performing regularization strategy by a relaxation of the 
costumer's demand leads to the reduction of the volatility of the 
sales turnover. 
 
\subsubsection{VaR approach on the dual thermal problem} 
 
In this subsection, we intend to exploit a stochastic model of the 
unavailability of the thermal plants in order to formulate a 
Value-At-Risk problem on the costs of thermoelectric power 
generation. Let $\ell$ be a thermal unit with $n_\ell$ thermal 
groups. Let $\alpha_{j,\ell}(t)$ be the probability that group $j$ 
of unit $\ell$ works at time step $t$ and $U_{j,\ell}^t$ the 
random variable such that $U_{j,\ell}^t=1$ if group $j$ works at 
time step $t$ and $U_{j,\ell}^t=0$ else. We suppose that the 
groups are regularly checked and, if necessary, repaired every $m_0$ 
time steps. Between two consecutive checking dates, we assume that 
the availability of the units is not changing. This means that 
between two consecutive checking dates, a given group is either 
working or it is out of work during the whole period. If $t_0=1$ 
and $t_k=m_0 k$ for $k \in \mathbb{N}^*$, then the probabilities 
$\alpha_{j,\ell}(t)$, for $t=t_k ,\ldots,t_{k+1}-1$ are the same 
and we only need to evaluate $\alpha_{j,\ell}(t_k), k \geq 0$. 
Those probabilities $\alpha_{j,\ell}(t_k)$ that a group $j$ of 
unit $\ell$ works at time step $m_0 k$ will depend on the past 
evolution of the availability of this group. If at time step 
$t_{k-1}$, the group was out of work, there is a big probability 
(say $1-\beta_1^\ell$ with $\beta_1^\ell$ small) that it works at 
time step $t_k$ (the time between two checking dates is greater 
than the mean time to repair) and a small probability 
$\beta_1^\ell$ that it is still out of work at time step $t_k$. 
Now if the group was working for the last $m$ periods delimited by 
the last $m+1$ checking dates, we can assume that the longer it 
has been working without failure (the larger $m$) the more likely 
it can break down at time step $t_k$. Thus, there is a decreasing 
function of $m$, $\beta_2^\ell(m)$ such that for any group j of unit 
$\ell$, 
$$ 
P(U_{j,\ell}^{t_k}=1|\mbox{Group j was working from} \;  t_{k-m}\; 
\mbox{to} \; t_{k-1})=\beta_2^\ell(m). 
$$ 
A particular case is the case where the state process of a given 
group is an homogeneous Markov chain where the state space is \{F,W\} 
where F stands for the failure state and W for the working state. 
In this case, $\beta_2^\ell(m)=\beta_2^\ell$ is fixed and 
corresponds to the probability for a group of unit $\ell$ to work 
on a given period knowing that it was working the period before. 
The transition matrix for the groups of unit $\ell$ is given by: 
$$ 
P_{\ell}=\left(\begin{array}{cc}\beta_1^\ell& 1-\beta_1^\ell \\ 
1-\beta_2^\ell&\beta_2^\ell \end{array}\right). 
$$ 
The probability $\alpha_{j,\ell}(t_k)$ is then given for $k \geq 
1$ by: 
$$ 
\alpha_{j,\ell}(t_k)=p_{F}^\ell(j)P_{\ell}^k(1,2) + 
p_{W}^\ell(j)P_{\ell}^k(2,2); 
$$ 
where $p_{W}^\ell(j)=1-p_{F}^\ell(j)=$ and $p_{W}^\ell(j)$ is the 
probability that group $j$ of unit $\ell$ works at the first time 
step. For the simplicity of the exposure, we assume that for a 
given unit $\ell$, either all the groups are working or all the 
groups are out of work at the first time step. Thus, 
$\alpha_{j,\ell}({t_k})$ is $j$-independent and 
$\alpha_{\ell}({t_k})$ will denote the probability that a group of 
unit $\ell$ works at time $t_k$. Further, we can partition the 
scenario tree in subtrees such that the root node and the leaves 
nodes of a given subtree respectively correspond to time steps 
$t_k$ and $t_{k+1}$ for some $k \in \mathbb{N}$. Thus the 
unavailability rates at the different nodes of any subtree of this 
partition are the same for a given unit. Let 
$\mathcal{O}=\displaystyle \cup_{k=1}^m \; \mathcal{O}_k$ be such that 
$\mathcal{O}_k$ are the nodes of the k-th subtree $\mathcal{S}_k$ 
in this partition. Let $T_k=\{(j,p),\;\;|\;\;j \; \in 
\mathcal{O}_k, p \; \in \mathcal{P}_j\}$. The dual thermal 
subproblem then writes: 
$$ 
\left\{\begin{array}{l} \min \;\; \displaystyle{\sum_{\ell}} 
\displaystyle{\sum_{k=1}^m} \;(c_{1 
\ell k} - \lambda_k)^T u_{t \ell k}\\ 
0 \leq u_{t \ell k} \leq \tau_{\ell}(k) \tau_{T}^\ell(k) 
P_{\max}^\ell d_k, 
\end{array} 
\right. 
$$ 
where $P_{\max}^\ell$ is the maximal available power of thermal 
unit $\ell$, $u_{t \ell k}=(u_{t \ell j p})_{\footnotesize{(j,p) 
\; \in \; T_k}}$, $c_{1 \ell k}=(c_{1 \ell j p 
})_{\footnotesize{(j,p) \; \in \; T_k}}$, 
$\lambda_k=(\lambda_{j,p})_{\footnotesize{(j,p) \; \in \; T_k}}$, 
$d_k=(d_{j,p})_{\footnotesize{(j,p) \; \in \; T_k}}$ (see the 
appendix), $\tau_{T}^\ell(k)$ gives the programmed unavailability 
rates for unit $\ell$ and the times subdivision of the set $T_k$ 
and $\tau_{\ell}(k)$ is the unavailability rate of unit $\ell$ for 
the nodes of the set $\mathcal{O}_k$. We reformulate this problem 
as the following problem: 
$$ 
\left\{\begin{array}{l} \min \;\; \displaystyle{\sum_{\ell}} 
\displaystyle{\sum_{k=1}^m} \tau_{\ell}(k) \;(c_{1 
\ell k} - \lambda_k)^T \tilde u_{t \ell k}\\ 
0 \leq \tilde u_{t \ell k} \leq \tau_{T}^\ell(k) P_{\max}^\ell 
d_k, 
\end{array} 
\right. 
$$ 
by setting $\tilde u_{t \ell k}=\frac{u_{t \ell 
k}}{\tau_{\ell}(k)}$. Notice that for quite a number of linear 
stochastic optimization problems, the random is only in the right 
hand side of the constraints as it is the case for the dual 
thermal subproblem. The above simple transformation allows to 
transfer the random in the objective and to implement a VaR method 
to compute robust solutions as was described in subsection 
\ref{var}. Given a confidence level $0 < \varepsilon <1$, we thus 
now introduce a VaR approach on the thermal plant cost/revenue 
balance: 
\begin{equation} \label{varthermal} 
(VaR)_{Benef}\left\{ 
\begin{array}{l} 
\min \;\; \gamma\\ 
P(\displaystyle{\sum_{\ell}} \displaystyle{\sum_{k}} 
\tau_{\ell}(k) \;(c_{1 
\ell k} - \lambda_k)^T u_{t \ell k} \leq \gamma) \geq 1 - \varepsilon,\\ 
0 \leq u_{t \ell k} \leq \tau_{T}^\ell(k) P_{\max}^\ell d_k. 
\end{array} 
\right. 
\end{equation} 
This problem may be understood as a problem of  {\bf maximization 
of the benefits} or equivalently as a problem of {\bf minimization 
of the losses}. We now need to study the modelling of the 
unavailability 
rates $\tau_{\ell}(k)$ to give an explicit form for problem \eqref{varthermal}.\\\\ 
\textbf{Modelling of the unavailability rates $\tau_{\ell}(k)$.} 
Let $\bar{P}^{\ell}_{\max}$ be the maximal power of a group in 
unit $\ell$. Then the theoretical maximal power available on the 
thermal unit $\ell$ is given by $P^{\ell}_{\max}=n_\ell 
\bar{P}^{\ell}_{\max}$. The maximal power available of unit $\ell$ 
for the nodes of the set $\mathcal{O}_k$ is then: 
\[ \tilde{P}^{\ell,k}_{\max}=\sum_{j=1}^{n_\ell} 
U_{j,\ell}^{t_k} \bar{P}^{\ell}_{\max}=n_\ell 
\bar{P}^{\ell}_{\max} {\sum_{j=1}^{n_\ell} U_{j,\ell}^{t_k} \over 
n_\ell}=P_{\max}^\ell\tau_{\ell}(k). 
\] 
If $t(k)$ is the time step associated with the root of the subtree 
$\mathcal{S}_k$, notice that under the above hypothesis, the 
random variable $n_\ell \tau_{\ell}(k)$ follows the binomial law $ 
{\cal B}(n_\ell,\alpha_{\ell}(t(k)))$. We then have 
$E[\tau_\ell(k)]=\alpha_{\ell}(t(k))$ and the variance of 
$\tau_\ell(k)$, 
$\mbox{Var}(\tau_\ell(k))={\alpha_{\ell}(t(k))(1-\alpha_{\ell}(t(k))) 
\over n_\ell}$. Now let $X_{\ell,k}$ be the random variable 
$\displaystyle{\sum_{\ell,k}} \; \tau_{\ell}(k) \;(c_{1 \ell k} - 
\lambda_k)^T u_{t \ell k}$. As the $(\tau_{\ell}(k))_{\ell,k}$ are 
independent we have: 
$$ 
E[X_{\ell,k}]=\displaystyle{\sum_{\ell,k}} \; 
\alpha_{\ell}(t(k))(c_{1 \ell k} - \lambda_k)^T u_{t \ell 
k},\;\;\mbox{and} \;\; 
VaR[X_{\ell,k}]=\displaystyle{\sum_{\ell,k}} \; u_{t \ell k}^T \, 
Q_{\ell k} \, u_{t \ell k}, 
$$ 
where the matrix $Q_{\ell k}$ is defined by 
$\frac{\alpha_{\ell}(t(k))(1-\alpha_{\ell}(t(k)))}{n_\ell}(c_{1 
\ell k} - \lambda_k)(c_{1 \ell k} - \lambda_k)^T$. From subsection 
\ref{var}, \eqref{varthermal} amounts to solve: 
\begin{equation}\label{dualtherm} 
(VaR)_{Benef}\left\{ 
\begin{array}{l} 
\min \;\;\displaystyle{\sum_{\ell,k}}  \alpha_{\ell}(t(k)) \;(c_{1 
\ell k} - \lambda_k)^T u_{t \ell k} + \kappa(\varepsilon)\sqrt{\displaystyle{\sum_{\ell,k}} \;u_{t \ell k}^T Q_{\ell k} u_{t \ell k}}\\ 
0 \leq u_{t \ell k} \leq \tau_{T}^\ell(k) P_{\max}^\ell d_k, 
\end{array} 
\right. 
\end{equation} 
where 
$\kappa(\varepsilon)=\sqrt{\frac{1-\varepsilon}{\varepsilon}}$. 
Notice that this is a second order cone optimization problem whose 
dimension will be high in practice (the number of nodes of the 
tree times the number of subdivision times). Another more 
conservative approach would be to use a VaR approach for each time 
subdivision of the nodes. This is possible as the dual thermal 
subproblem is separable with respect to the time subdivisions of 
the nodes. The advantage of this approach is that we have an 
explicit solution for the new dual thermal subproblem and the new 
dual thermal subproblem is again separable with respect to both 
the thermal units and the time subdivisions. This allows to solve 
problems of big sizes. Moreover, the unavailability rate of unit 
$\ell$ for each time subdivision follows a binomial law which can 
be approximated by a Gaussian law if ($n_\ell \, \alpha_\ell(t(k)) 
\geq 10$ and $n_\ell \, (1-\alpha_\ell(t(k))) \geq 10)$ or 
$n_\ell$ big enough, say $n_\ell \geq 6$. At last, we mention that 
a particular case would be to consider that $\alpha_\ell(t_k)$ 
doesn't depend on $k$. In what follows we both suppose that 
$\alpha_\ell(t_k)$ doesn't depend on $k$ and that the VaR approach 
on the thermal subproblem is done at each time subdivision. This 
will allow to check how the VaR approach on the thermal 
subproblems gives an immunization with respect to the uncertainty 
we have on the 
unavailability of the units which works surprisingly very well in practise.\\\\ 
{\bf Economical Interpretation.} Qualitatively, one can say that the 
objective of ($VaR_{Benef}$) is to maximize the benefits of a 
thermal unit while ensuring some robustness with respect to plants 
random unavailability. 
 
\subsubsection{Intermediate summary} 
 
At this point, we have introduced two different regularizations on 
the original dual optimization problem that was initially 
formulated as : 
\begin{equation} 
(DP) \qquad \displaystyle \max_{\lambda} {\theta}_d(\lambda) + \max_{\lambda}\left ({\theta}_T+ {\theta}_N + \theta_H + \theta_ J  )\right.(\lambda) 
\end{equation} 
where  the  $\theta_{j}$ are given by \eqref{dualdec}. 
\begin{enumerate} 
\item From a relaxation of the demand, we formulate the dual regularized problem 
\begin{equation} 
(VaR_{FA}) \qquad \displaystyle \max_{\lambda} 
{\theta}_d^{\RC}(\lambda) + \max_{\lambda}\left ({\theta}_T+ 
{\theta}_N + \theta_H + \theta_ J  )\right.(\lambda) 
\end{equation} 
 with 
\begin{equation} 
\begin{array}{lll} 
{\theta}_d^{\RC}(\lambda)&=& E_{\omega}(d(\omega))^T \lambda - 
\kappa(\varepsilon_1) \sqrt {\lambda^T \Gamma \lambda}. 
\end{array} 
\end{equation} 
where $\varepsilon_1$ is the confidence level chosen to implement 
$VaR_{FA}$. 
\item From a rewriting of the random unavailability of the thermal plants, we 
get 
\begin{equation} 
(VaR_{Benef}) \qquad \displaystyle \max_{\lambda}  {\theta}_d(\lambda)  + \max_{\lambda} {\theta}_T^{\RC}(\lambda) + \max_{\lambda}\left ( {\theta}_N + \theta_H + \theta_ J  )\right.(\lambda)  \end{equation} 
 with 
\begin{equation} \label{versionvartherm2} 
\left\{ \begin{array}{l} 
{\theta}_T^{\RC} (\lambda)  = \displaystyle{\sum_{n,p,\ell \in \mathcal{L}_T}}  \min_{0 \leq u \leq \tau_{T,n}^\ell P_{\max}^\ell d_{n,p}} V_\ell(u,\lambda_{n,p}) \\ 
V_\ell(u,\lambda_{n,p})=\left(\alpha_\ell (c_{1 \ell n 
p}-\lambda_{n,p}) + \kappa(\varepsilon_2) \sqrt{\frac{\alpha_\ell 
(1- \alpha_\ell)}{n_\ell}} |c_{1 \ell n p}- 
\lambda_{n,p}|\right)u. 
\end{array} 
\right. 
\end{equation} 
\end{enumerate} 
where $\varepsilon_2$ is the confidence level chosen to implement 
$VaR_{Benef}$. If the we use a VaR technique on the whole thermal 
subproblem then ${\theta}_T^{\RC} (\lambda)$ is given by 
\eqref{dualtherm}. Finally, we can also combine the two previous 
regularizations as the mixt problem: 
\begin{equation} 
(VaR_{mixt}) \qquad \displaystyle{\max_{\lambda}} \;\; 
{\theta}_{\mbox{mixt}}^{\RC} (\lambda)=\displaystyle 
\max_{\lambda} {\theta}_d^{\RC}(\lambda)  + \max_{\lambda} 
{\theta}_T^{\RC}(\lambda) + \max_{\lambda}\left ( {\theta}_N + 
\theta_H + \theta_ J )\right.(\lambda). \end{equation}

\section{Implementation and Numerical simulations} \label{testsEDF} 
 
\subsection{Implementation and simulation protocol} 
 
To solve the primal optimization problem \eqref{base3}, we have to 
solve \begin{math}\displaystyle \min_{u\in 
\UC}\max_{\lambda}L(u,x_h,x_e,\lambda) \end{math}, where $L$ is the 
usual Lagrangian. This  will be equivalent to the dual problem 
$\displaystyle{\max_{\lambda} \, \theta(\lambda)}$  if only 
thermal and hydro units are considered. Indeed, in this case, 
problem \eqref{base3} is a below bounded linear program and both 
the primal and the dual are equivalent to each other. If we take 
into account EJP contracts, the set of constraints is not convex 
and the duality gap is strictly positive. However, the weak 
duality relationship still holds: 
\begin{equation} 
 \displaystyle \min_{u\in \UC}\max_{\lambda}L(u,x_h,x_e,\lambda) \geq \displaystyle{\max_{\lambda} \; \theta(\lambda)}. 
\end{equation} 
Moreover, numerical simulations have shown that the duality gap is 
generally quite small. The dual problem thus allows us to 
approximate primal solutions and estimate marginal prices. 
 
\subsubsection{Space decomposition for Optimization} \label{spacedec} 
 
First, we describe the space decomposition method for  the dual 
function $\theta$ and next, we explain the adaptations necessary 
for the regularized problems. The dual function $\theta$ is non 
differentiable, concave and separable with respect to the units as 
it writes $\displaystyle{\theta(\lambda)=\lambda^T d + \sum_{\ell 
\in \mathcal{L}} \; \theta^\ell(\lambda)}$ with 
$$ 
\displaystyle{\theta^\ell(\lambda)=\min_{(x_{\bullet}^\ell(\bullet),u_{\bullet}^\ell(\bullet)) 
\in \chi_\ell} \; \sum_{n \in \mathcal{O}} \,\sum_{p \in 
\mathcal{P}_n}} \, \pi_n \, 
{\mathcal{C}}_{n,p}^\ell(x_{n}^\ell(p),u_{n}^\ell(p))-\lambda_{n,p} 
P_{n,p}^{\ell}(x_n^\ell(p),u_n^\ell(p)), 
$$ 
the dual function of the subproblem associated with unit $\ell$. 
This is especially of interest to treat problems of big size as it 
is the case for our application. To maximize $\theta$ (or which is 
the same to minimize the convex function $-\theta$) we use a 
bundle method described in \cite{claude12}. This requires to build 
a black box which, for any $\lambda \in \mathbb{R}^D$ is able to 
compute $-\theta(\lambda)$ and to give an arbitrary subgradient 
$s(\lambda) \in 
\partial(\footnotesize{-}\theta(\lambda))$. Let $L_{\ell}$ be the partial 
Lagrangian associated with the partial dual function corresponding 
to unit $\ell$. The computation of $-\theta(\lambda)$ is done 
solving the different optimization problems associated with the 
different production units. As for a computation of a subgradient, 
if $y^\ell(\lambda)=\mbox{argmin}_{x,u} \; 
L_\ell(x_{\bullet}^\ell(\bullet),u_{\bullet}^\ell(\bullet),\lambda)$ 
and 
$P^\ell(\lambda)={(P^\ell_{n,p}(y_{n,p}^\ell(\lambda)))}_{n,p}$ 
then 
$$ 
-\theta^\ell(\mu) \geq -\theta^\ell(\lambda) + {(\mu-\lambda)}^T 
\,P^\ell(\lambda) \;\; \textrm{for all} \;\;(\lambda,\mu) \in 
\mathbb{R}^D, 
$$ 
which shows that $P^\ell(\lambda) \in 
\partial(\footnotesize{-}\theta^\ell(\lambda))$. We then 
immediately have 
$$ 
s(\lambda)=-d+ \displaystyle{\sum_{\ell \in \mathcal{L}} \; 
P^\ell(\lambda)} \in \partial(\footnotesize{-}\theta(\lambda)). 
$$ 
More precisely, following \cite{bgls}, a global resolution by an 
iterative scheme can be described with 4 steps starting from a 
reference price $\lambda$ used to initialize the algorithm with 
index  $k=1$ and $\lambda_1 =\lambda$. 
\begin{enumerate} 
\item At  iteration $k$, decomposition in subproblems and computation of the  local solution of subproblem  $\ell$: $y^\ell(\lambda_k)$; 
\item Evaluation of the dual function  $\theta$   at the point  $\lambda_k$ and computation of a subgradient $s(\lambda_k)$; 
\item Updating of the multipliers by the  coordinator using a black box method (i.e computation of $\lambda_{k+1}$); 
\item Updating of the index  $k \leftarrow k+1$ and go to step 1. 
\end{enumerate} 
The robustifications proposed still remain within the same 
framework, so we can solve them in a space decomposition 
framework. The adaptation of the above algorithm to the 
robustifications proposed is easy since we just need to modify 
step 2. 
\begin{itemize} 
\item For the dual regularized problem $(VaR_{FA})$, we have to increase the value of the dual function 
by  $\displaystyle  -\kappa(\varepsilon_1) \sqrt{ {\lambda}^T 
\Gamma \lambda }$ and the sugradient by $\displaystyle 
-\kappa(\varepsilon_1) \frac{\Gamma \lambda}{{\lambda}^T \Gamma 
\lambda}$ where $\varepsilon_1$ is the confidence level used to 
implement $(VaR_{FA})$. So the only extra cost of this model is 
the estimation of matrix $\Gamma$. 
\item The method $VaR_{benef}$ simply modifies the thermal dual problem which becomes 
problem \eqref{varthermal} (if the VaR approach is done on the 
whole tree) or \eqref{versionvartherm2} else. If 
\eqref{versionvartherm2} is used, then the solution of the new 
thermal dual problem is still separable with respect to both the 
units $\ell$ and the time subdivisions $p$ and the optimal 
commands $u_{n,p}^\ell$ for time subdivision $p$, node $n$, unit 
$\ell$ are immediately given by the following formulas: 
\begin{equation} 
\left\{ 
\begin{array}{ll} 
u_{n,p}^\ell=0 & \mbox{ if } \left(\pi_n c_\ell \geq \lambda_{n,p} 
\mbox{ or } \pi_n c_\ell < \lambda_{n,p} 
\mbox{ and } \displaystyle {\alpha_\ell \leq \frac{{\kappa(\varepsilon)}^2}{{\kappa(\varepsilon)}^2 + n_\ell}}\right),\\ 
u_{n,p}^\ell=\tau_{T,n}^\ell \, P_{\max}^\ell \, d_{n,p} & \mbox{ 
else.} 
\end{array} 
\right. 
\end{equation} 
\end{itemize} 
 
\subsubsection{Data and simulation protocol} 
 
The data used for the simulations are inspired from real data. We 
suppose that each time step is divided in $L$ time subdivisions 
also called hourly posts. Our generation strategy will be tested 
on a set of 456 scenarios. To each scenario is associated a 
realization of the inflows in the reservoirs, the unavailability 
rates of the thermal units and of the demand at each time 
subdivision of the year. From these scenarios we build 3 different 
trees. Each tree corresponds to a vision more or less difficult of 
the evolution of the inflows and the demand on the coming year. We 
will call those trees Easy tree, Median tree and Difficult tree 
with evident interpretation of the predictions of the demands and 
the inflows on those  trees. The scenario trees are trees of depth 
364 days, with 5227 nodes. At each node, we know the demand vector 
of the demands for all the posts $p$ of this node, the inflows for 
all the hydro reservoirs and the programmed unavailability rates 
of all the hydro and thermal units. There are $L=3$ hourly posts 
per day. We use the following generation units: 
\begin{itemize} 
\item Eleven thermal units. Every thermal unit $\ell$ is described by its (unitary) generation cost, its maximal and minimal power, the number of thermal groups and the probability $\alpha_\ell$ that a group works. 
\item Two independent hydro plants. Each hydro plant is connected to a different reservoir. We know the maximal stock (in GWh) of each reservoir, the initial stock of each reservoir and the maximal power (in MW) of each plant. The maximal stock of the biggest reservoir is around 30 times that of the other reservoir. This explains why we will essentially be interested in the evolution of the biggest reservoir stock on the year. 
\item An EJP contract of 22 days with maximal available power: $P_{J}^1=2467$ MW. 
\end{itemize} 
Before presenting the results it remains to explain how the covariance matrix $Q$ involved in $VaR_{FA}$ method is chosen. 
\subsubsection{Calibration of the covariance matrix $Q$} 
Since all we had were the scenarios of demands and the demands at the nodes of the three different trees generated from those scenarios following the lines of \cite{carpentier1}, it was difficult to calibrate the matrix $Q$. However, to have an idea of the impact of this method on the simulation process we supposed the demand at the different hourly posts were uncorrelated. We thus had to deem the diagonal elements of $Q$ corresponding to the variances $\sigma^2(n,p)$ of the demand for node $n$ and post $p$. This node $n$ is associated to a time step $t$. To estimate $\sigma(n,p)$ we sort the demands of this time step by increasing order. Let's denote by $d_{(i)}^t, 1\leq i \leq m$ the sample of demands ordered by increasing order for time step $t$. We then choose for $\sigma(i)$ which determines the uncertainty we have on $d_{(i)}^t$: 
$$\sigma(i)=\min(\frac{d_{(i+1)}^t 
-d_{(i)}^t}{2},\frac{d_{(i)}^t -d_{(i-1)}^t}{2}), 
$$ 
with the convention $d_{(0)}^t=0$ and $d_{(m+1)}^t=2 d_{m}^t - d_{(m-1)}^t$. We will compare the results we obtained with the nominal model and with the VaR approaches $VaR_{FA}$ and 
$VaR_{benef}$. The outputs we are interested in are guided by the defaults of the nominal model outlined in subsection \ref{default}. We are thus interested in the distribution of the simulated costs and in the behavior of the hydro reservoirs. 
\subsection{Numerical results} 
\subsubsection{Central and dispersion characteristics of the costs.} 
We provide the mean and the standard deviation of the simulated costs. We also give the empirical quantile of order 0.95 (VaR 5\%) and of order 0.99 (VaR 1\%) of the distribution of these costs. We give the results using the three different trees (Easy, Median and Difficult tree) to solve the optimization problem and for all the methods. We also take a look at the method (called Mixt) consisting in cumulating the two modifications proposed by the two VaR approaches. 
 
\begin{center} 
\begin{tabular}{|c|c|c|c|c|} 
\hline 
Output & Nominal  & Robust $VaR_{FA}$ & Robust $VaR_{benef}$ & Mixt\\ 
\hline 
\begin{tabular}{c} 
Mean\\ 
St. Dev. 
\end{tabular} 
& 
\begin{tabular}{r} 
467 529 291\\ 
47 864 238 
\end{tabular} 
& 
\begin{tabular}{r} 
488 271 561\\ 
46 609 492 
\end{tabular} 
& 
\begin{tabular}{r} 
459 515 733\\ 
31 750 330 
\end{tabular} 
& 
\begin{tabular}{r} 
459 991 109\\ 
30 597 481 
\end{tabular}\\ 
\hline 
\begin{tabular}{r} 
VaR 1\%\\ 
VaR 5\% 
\end{tabular} 
& 
\begin{tabular}{r} 
671 599 214\\ 
543 786 128 
\end{tabular} 
& 
\begin{tabular}{r} 
672 311 733\\ 
574 245 712 
\end{tabular} 
& 
\begin{tabular}{r} 
557 798 070\\ 
517 836 912 
\end{tabular} 
& 
\begin{tabular}{r} 
558 719 609\\ 
518 856 055 
\end{tabular}\\ 
\hline 
\end{tabular} 
\end{center} 
 
\begin{center} 
\textbf{Central and dispersion characteristics of the empirical distribution of the simulated costs on the Easy tree (management horizon of 1 year).} 
\end{center}

\begin{center} 
\begin{tabular}{|c|c|c|c|c|} 
\hline 
Output & Nominal  & Robust $VaR_{FA}$ & Robust $VaR_{benef}$ & Mixt\\ 
\hline 
\begin{tabular}{c} 
Mean\\ 
St. Dev 
\end{tabular} 
& 
\begin{tabular}{r} 
466 498 435\\ 
46 540 234 
\end{tabular} 
& 
\begin{tabular}{r} 
486 587 235\\ 
47 075 591 
\end{tabular} 
& 
\begin{tabular}{r} 
462 184 762\\ 
29 154 061 
\end{tabular} 
& 
\begin{tabular}{r} 
462 334 138\\ 
28 857 445 
\end{tabular}\\ 
\hline 
\begin{tabular}{r} 
VaR 1\%\\ 
VaR 5\% 
\end{tabular} 
& 
\begin{tabular}{r} 
689 406 053\\ 
548 912 439 
\end{tabular} 
& 
\begin{tabular}{r} 
679 536 989\\ 
573 426 666 
\end{tabular} 
& 
\begin{tabular}{r} 
557 838 702\\ 
516 533 109 
\end{tabular} 
& 
\begin{tabular}{r} 
554 297 141\\ 
515 210 039 
\end{tabular}\\ 
\hline 
\end{tabular} 
\end{center} 
 
\begin{center} 
\textbf{Central and dispersion characteristics of the empirical distribution of the simulated costs on the Median tree (management horizon of 1 year).} 
\end{center} 
 
\begin{center} 
\begin{tabular}{|c|c|c|c|c|} 
\hline 
Output & Nominal  & Robust $VaR_{FA}$ & Robust $VaR_{benef}$ & Mixt\\ 
\hline 
\begin{tabular}{c} 
Mean\\ 
St. Dev. 
\end{tabular} 
& 
\begin{tabular}{r} 
464 712 495\\ 
44 615 708 
\end{tabular} 
& 
\begin{tabular}{r} 
479 773 359\\ 
47 172 746 
\end{tabular} 
& 
\begin{tabular}{r} 
465 886 909\\ 
28 146 821 
\end{tabular} 
& 
\begin{tabular}{r} 
465 097 476\\ 
28 272 953 
\end{tabular}\\ 
\hline 
\begin{tabular}{r} 
VaR 1\%\\ 
VaR 5\% 
\end{tabular} 
& 
\begin{tabular}{r} 
667 480 170\\ 
543 728 341 
\end{tabular} 
& 
\begin{tabular}{r} 
693 735 355\\ 
562 040 194 
\end{tabular} 
& 
\begin{tabular}{r} 
556 573 568\\ 
517 110 690 
\end{tabular} 
& 
\begin{tabular}{r} 
556 464 708\\ 
518 142 699 
\end{tabular}\\ 
\hline 
\end{tabular} 
\end{center} 
\begin{center} 
\textbf{ Central and dispersion characteristics of the empirical distribution of the simulated costs on the Difficult tree (management horizon of 1 year).} 
\end{center} 
 
We will essentially notice the following points: 
\begin{itemize} 
\item The average managing costs on the whole scenarios and for each of the three trees are quite close for the four methods considered. For method $VaR_{FA}$, the average costs are between 3.2\% and 4.4\% greater than the  nominal method. For method $VaR_{benef}$ and Mixt, the costs can be greater or less than those of the nominal method. Nevertheless, those costs only vary between 0.08\% and 1.7\% compared with those of the nominal model. 
\item Method $VaR_{FA}$ does not go into the good sense as the standard deviation of the costs as well as the VaR at 1\% and 5\% are greater than the same quantities computed for the nominal model. As for the methods $VaR_{benef}$ and Mixt they lead in all the cases to reductions of the standard deviation of the costs (till 38\% of reduction on the Difficult tree) and of the VaR at 1\% and 5\%. A reason that could explain the bad results of method $VaR_{FA}$ would be that the relaxation of the demand constraint in an ellipsoid works as an opportunity for the system to have an another reserve to perform its optimisation. But, this reserve does not exist in the Monte-Carlo simulation and thus, the strategy reveals itself to be too optimistic. The good results of $VaR_{benef}$ could be attributed to the fact that the thermal problem is the optimization problem really dimensioning (the total thermal costs are around 100 times greater than the maximal valorization possible of the biggest reservoir). It thus indeed seems 
interesting to envisage a robust approach which takes into account the only random involved in the thermal subproblem : the unavailability rates of the thermal plants. 
\end{itemize} 
 
 
\subsubsection{Study of the trajectories of the biggest reservoir.}

\begin{figure}[H] 
\begin{center} 
\begin{tabular}{ccc} 
\includegraphics[angle=0, width=4.5cm]{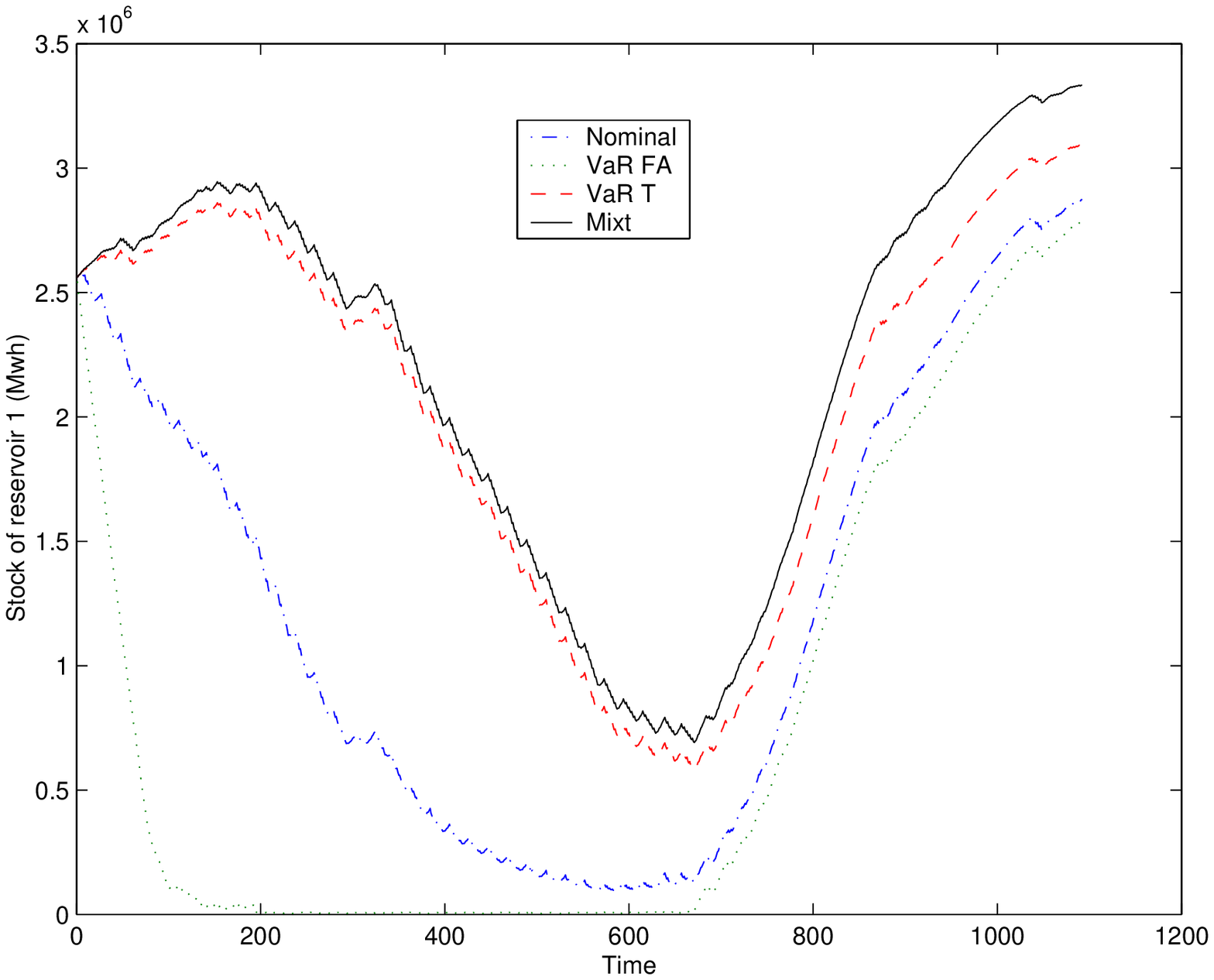} 
& 
\includegraphics[angle=0,width=4.5cm]{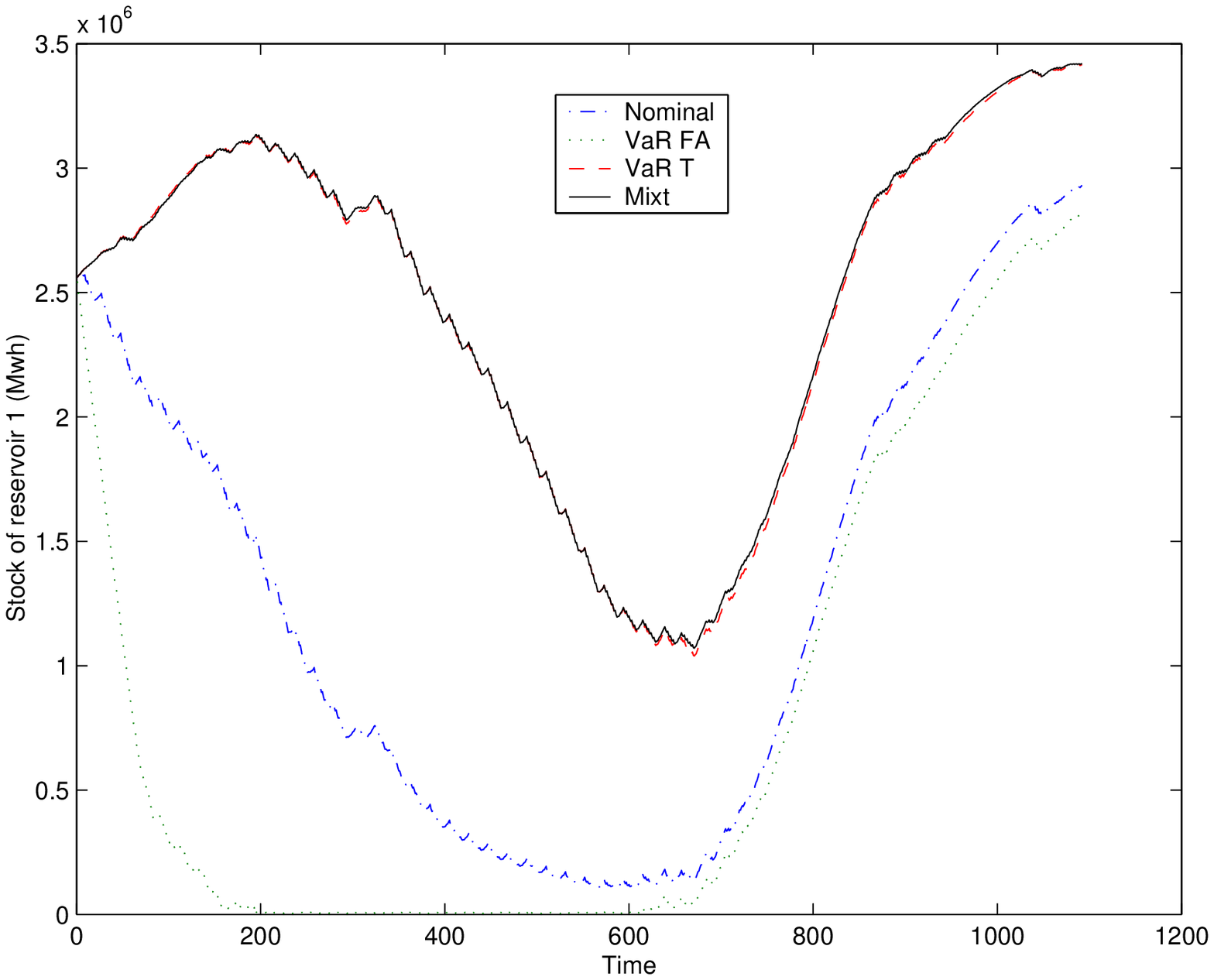} 
& 
\includegraphics[angle=0,width=4.5cm]{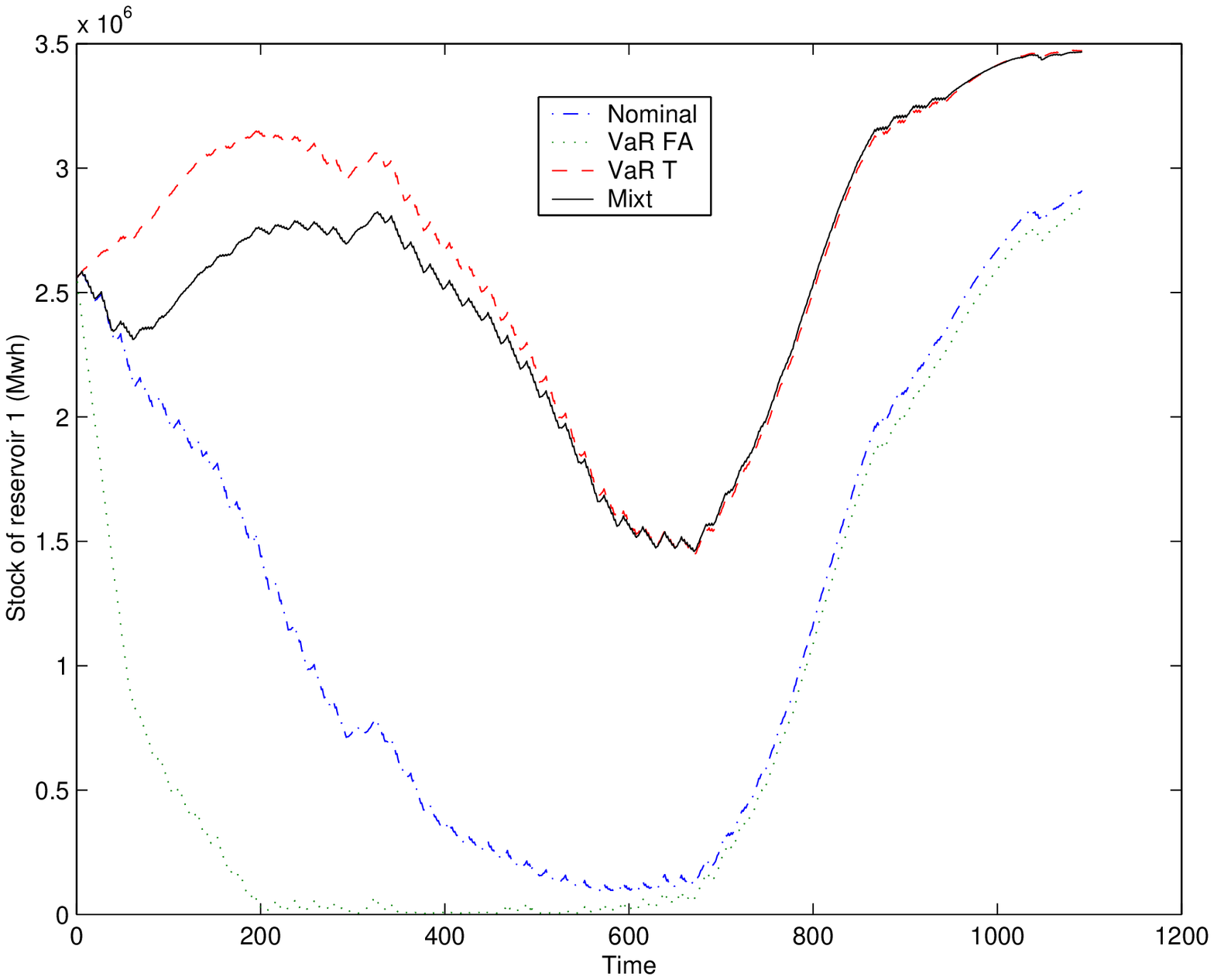}\\ 
Easy tree & Median tree& Difficult tree 
\end{tabular} 
\caption{Evolution of the average stock (in MWh and on the whole scenarios) of the biggest reservoir during the year for all the methods and using the Easy, Median and Difficult trees as a support of the optimization process.} 
\end{center} 
\end{figure} 
For a given method, the strategy of management of the biggest reservoir softly vary when we change the scenario tree. The nominal and $VaR_{FA}$ methods  tend to empty more the reservoir whose level increases at the end of the year. The methods $VaR_{benef}$ and Mixt do not use the reservoir or very little at the beginning of the year. Globally, the reservoir has a higher level with those methods and is nearly full (its maximal stock is 3500 GWh) at the end of the year. In what follows, we say that the biggest reservoir (res for short) is at a low level if it contains at most 5\% of its maximal stock. We say that a reservoir is at a high level if it attains a level greater or equal to its level at the beginning of the year less 5\% of its maximal stock. Using those notations the array below permits to precise tendencies already observed on the above curves. 
\begin{center} 
\begin{tabular}[b]{|c||cccc|cccc|} 
\hline 
{\small \# of weeks} & \multicolumn{4}{c|}{High level res} & \multicolumn{4}{c|}{Low level res}\\ 
& Init & $VaR_{FA}$ & $VaR_{benef}$ & Mixt &  Init & $VaR_{FA}$ & $VaR_{benef}$ & Mixt\\ 
\hline 
1&437  & 405  &456   &456        &426  &456  & 5   & 9\\ 
2&421  & 387  &456   &456      &423  &456  & 4   & 4\\ 
3&408  & 377  &456   &456      &423  &456  & 3   & 3\\ 
4&390  & 344  &456   &456      &417  &456  & 3   & 2\\ 
5&375  & 312  &456   &456      &412  &456  & 3   & 2\\ 
10&201 & 64  &456   &449      &365  &456  & 0   & 0\\ 
15&69  & 0  &449   &410      &256  &456    & 0   & 0\\ 
20&29  & 0  &438   &373       &92   &453    & 0   & 0\\ 
25&10  & 0  &414   &329      &14   &257    & 0   & 0\\ 
30&5   &  0  &371   &268      &0    &25    & 0   & 0\\ 
\hline 
\end{tabular} 
\end{center} 
 
\begin{center} 
Number of scenarios among 456 for which the biggest reservoir is at least $X$ weeks with a high or low level (optimizer launched on the Difficult tree and with all the methods). 
\end{center}

 
\subsubsection{Comparison of the distribution of the costs.} 
\begin{figure}[H] 
\begin{center} 
\begin{tabular}{cc} 
\includegraphics[angle=0, width=6cm]{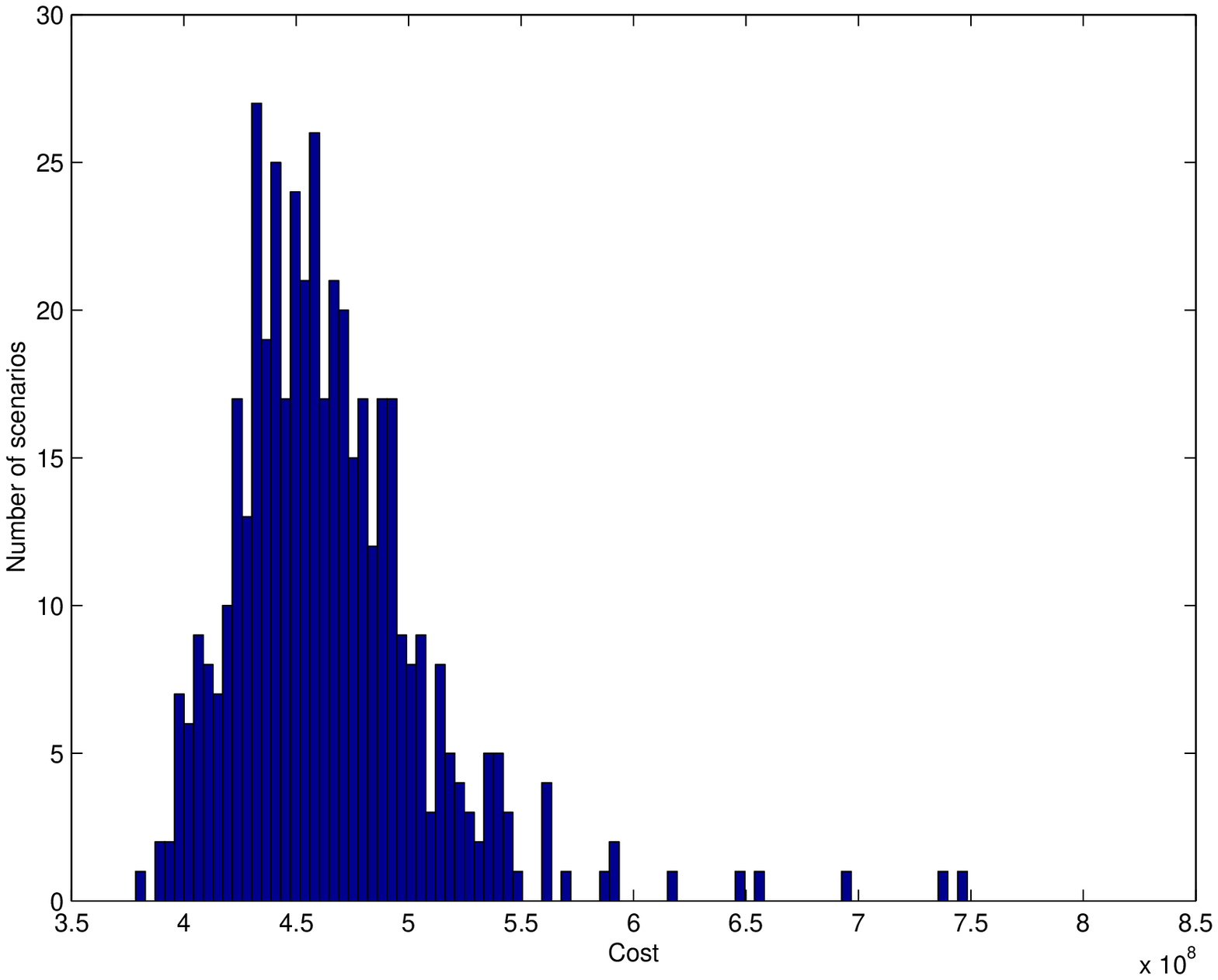} 
& 
\includegraphics[angle=0, width=6cm]{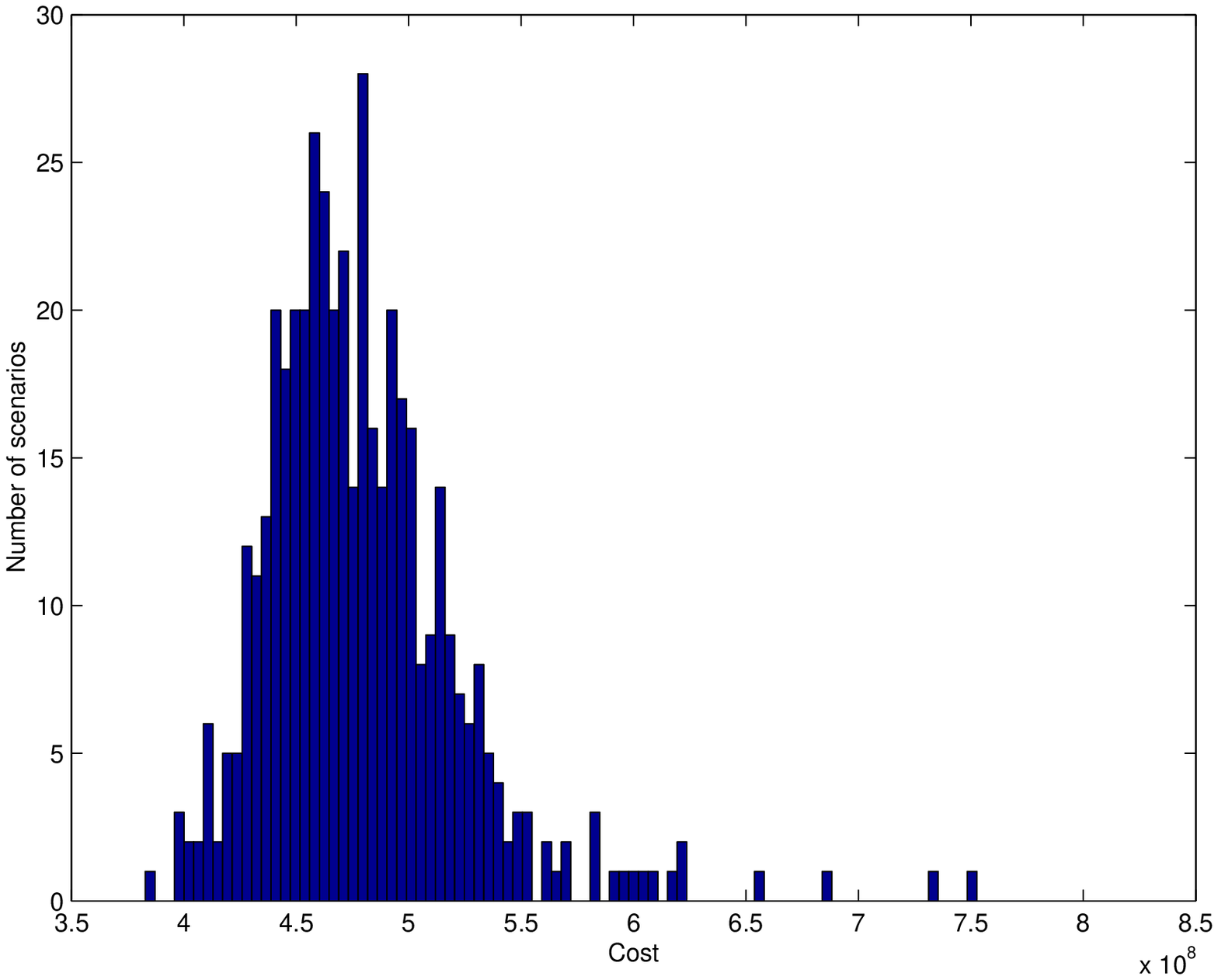}\\ 
Nominal & Robust $VaR_{FA}$ 
\end{tabular} 
\end{center} 
\end{figure} 
\begin{figure}[H] 
\begin{center} 
\begin{tabular}{cc} 
\includegraphics[angle=0, width=6cm]{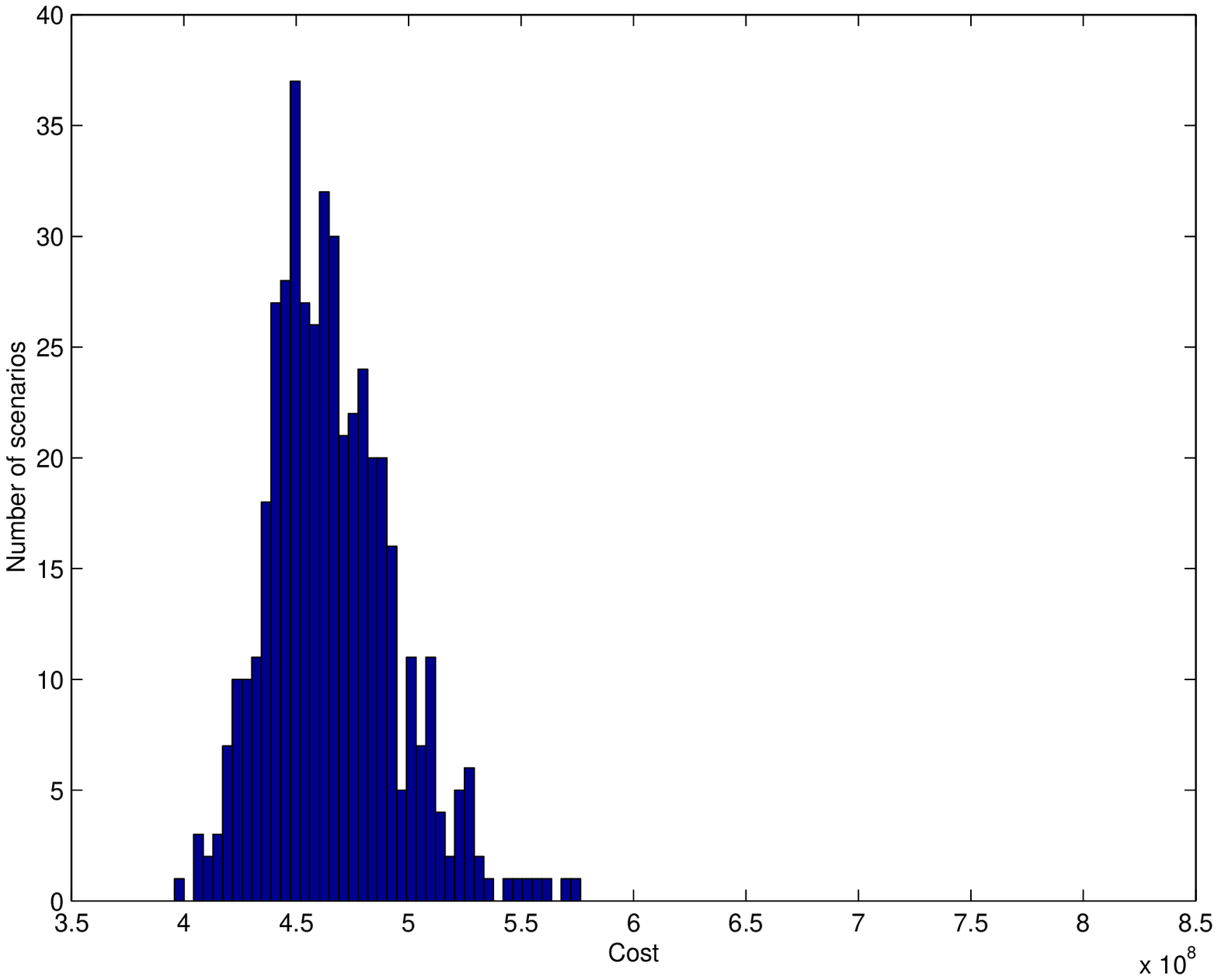} 
& 
\includegraphics[angle=0, width=6cm]{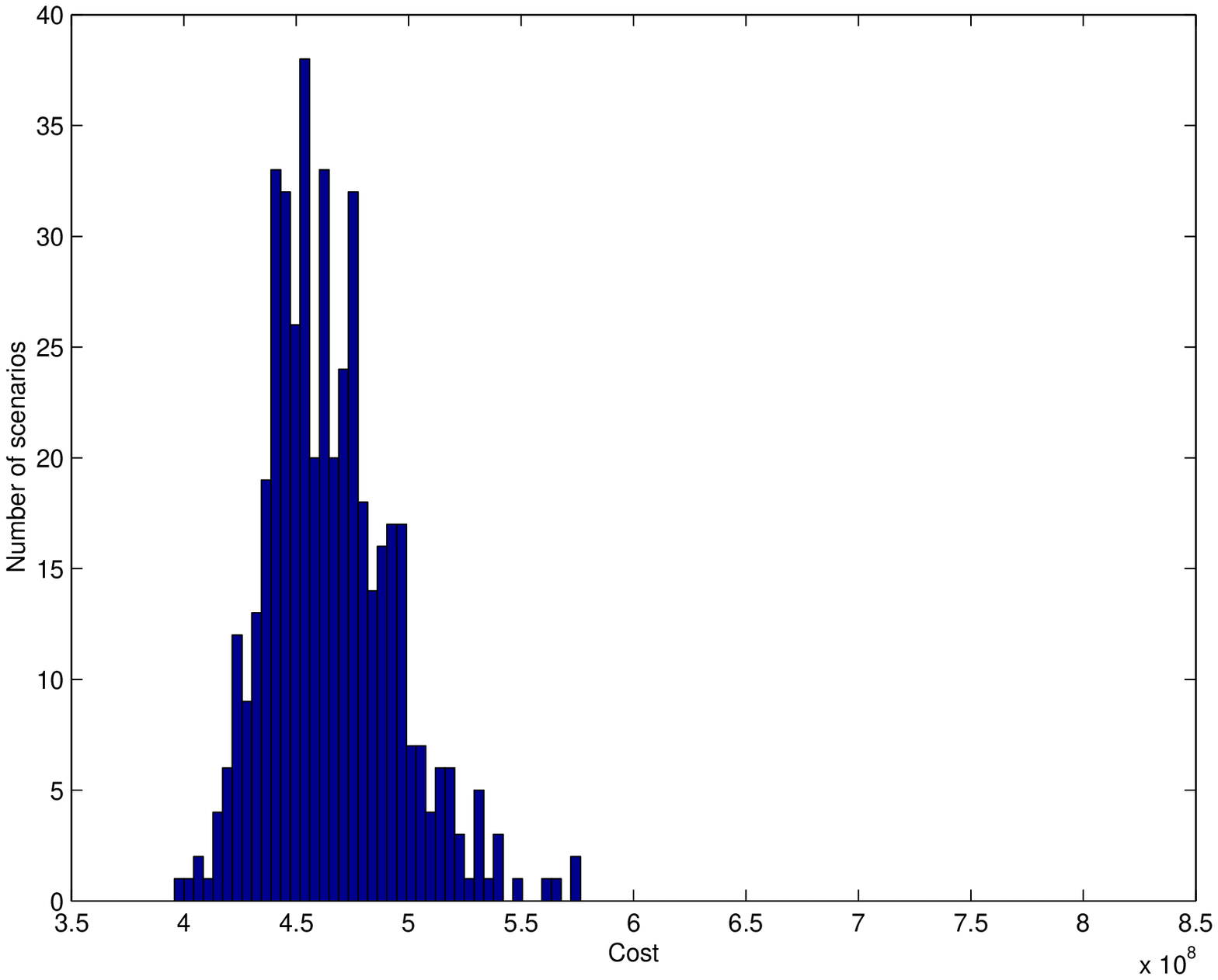}\\ 
Robust $VaR_{benef}$ & Robust Mixt 
\end{tabular} 
\end{center} 
\begin{center} 
\caption{Empirical densities of the management costs on the Difficult tree and for each method (management horizon 1 year).} 
\end{center} 
\end{figure} 
The empirical densities of the management costs for the nominal and  $VaR_{FA}$ methods have tails of distribution bigger than those of methods $VaR_{benef}$ and Mixt. A few scenarios are of very high cost for the nominal and $VaR_{FA}$ methods. On the contrary, the dispersion of the costs for the models $VaR_{benef}$ and Mixt is smaller.  We can illustrate those words by a few figures: 
\begin{itemize} 
\item The scenario of highest cost corresponds to costs of $7.448*10^8$,$8.232*10^8$,$5.738*10^8$ and $5.753*10^8$ for respectfully the nominal $VaR_{FA}$, $VaR_{benef}$ and Mixt methods. 
\item The less costly scenario for each of the models have close costs that are worth $3.822*10^8$,$3.868*10^8$,$3.997*10^8$ and $3.991*10^8$ for respectfully the nominal,$VaR_{FA}$,$VaR_{benef}$ and Mixt methods. 
\end{itemize} 
The shape below of the empirical cumulative distribution function of the management costs confirms this tendency. The nominal method has a bigger density of scenarios whose cost is less than $4.7*10^8$. On the other hand, there exists a non neglectable number of scenarios of high costs. For $VaR_{FA}$, it is worse as the density of scenarios whose cost is less than $4.7*10^8$ is small. 
\begin{figure}[H] 
\begin{center} 
\begin{tabular}{cc} 
\includegraphics[angle=0, width=6cm]{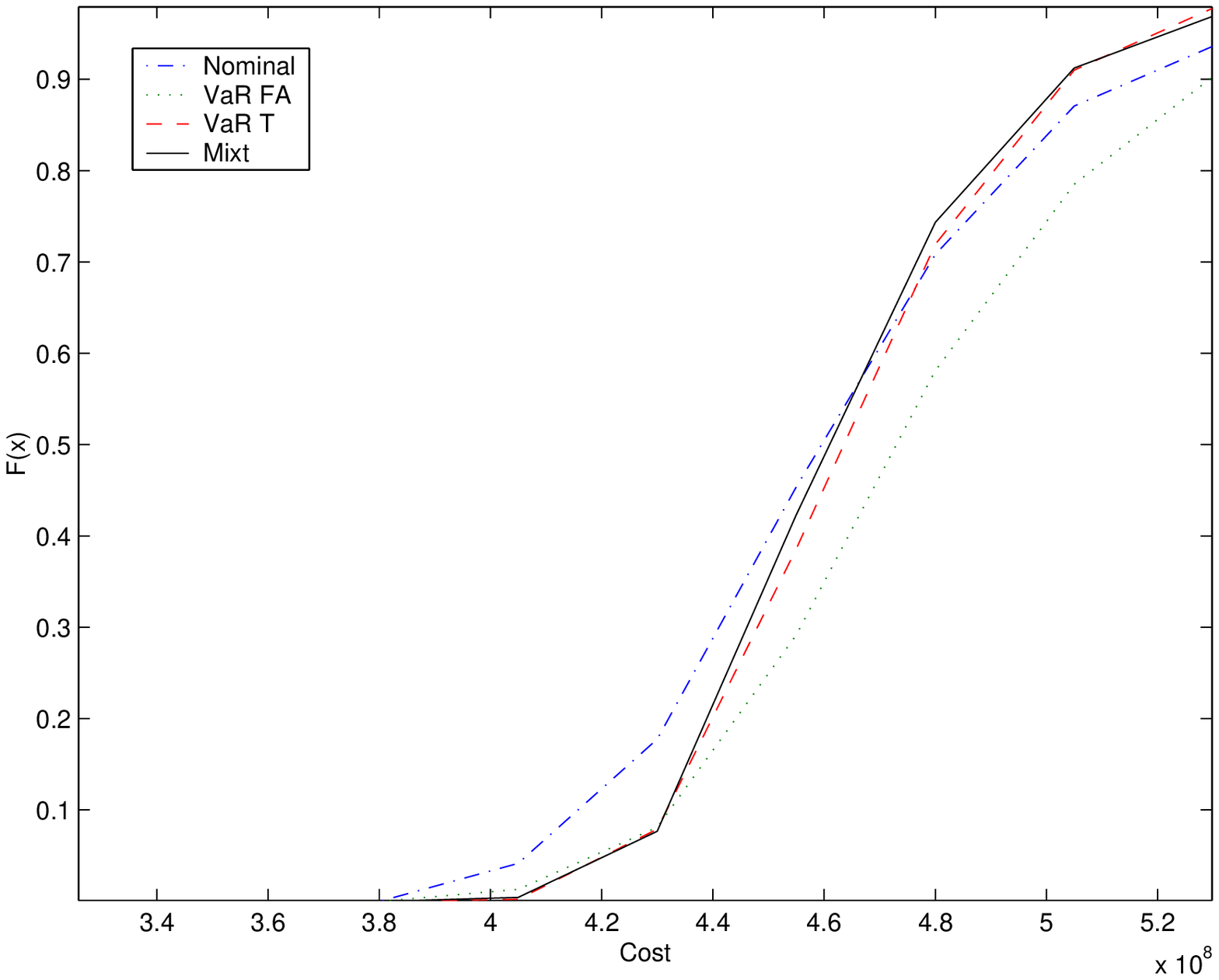} 
& 
\includegraphics[angle=0, width=6cm]{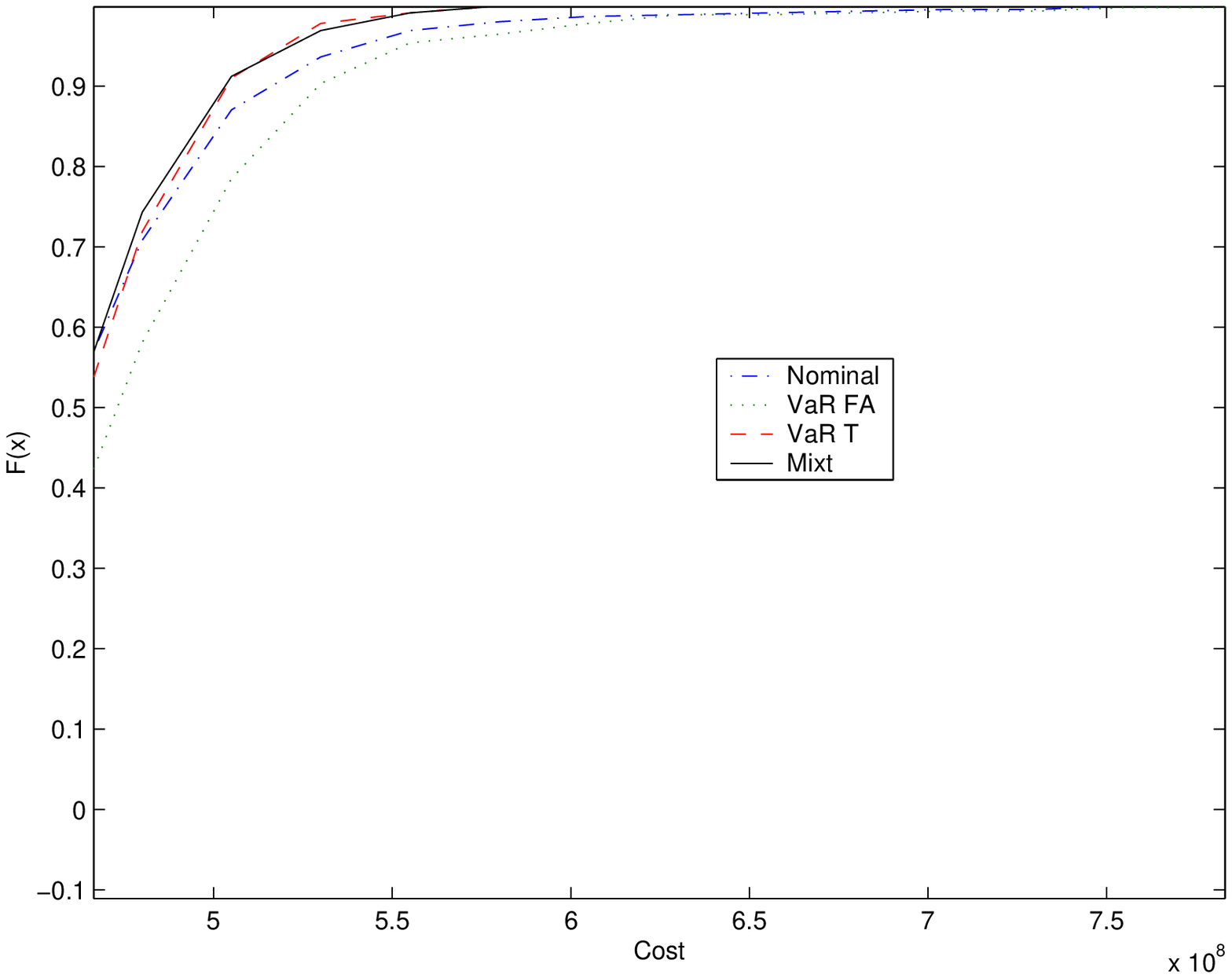}\\ 
Lower part & Higher part 
\end{tabular} 
\caption{Empirical cdf of the costs on the Difficult tree for the four methods.} 
\end{center} 
\end{figure} 
 
\section{Conclusion} 
 
This paper has presented an application of Value-At-Risk methods 
to robustify a stochastic optimization problem of yearly electric 
generation management. The starting point of our investigation was 
a nominal model which conducted to a big standard deviation of the 
costs and which emptied too much the hydro reservoirs. Two models 
have been proposed to reduce those drawbacks. If the first model 
has not been concluding in practice, the second model has revealed 
very adapted to the objectives. On the one hand, this model 
conducts to a diminution of the standard deviation of the 
simulated management costs. On the other hand it tends to empty 
less the biggest reservoir. The success of this model comes from 
the fact that it modifies the thermal problem which is the 
optimizing problem really dimensioning. The model $VaR_{FA}$, as 
for itself, is a classical robustification of the dual problem 
which takes into account the uncertainty on the demand. From a 
theoretical point of view, it permits to somehow stabilize the 
Lagrange multipliers which correspond to electricity prices in our 
application. On the other hand, this method should be interesting 
on difficult scenarios. To improve the results of this method, the 
problem of estimation of the matrix $\Gamma$ should be carefully 
studied. 
 
\paragraph{Acknowledgement}  
The authors are grateful to Anatoli Iuuditski of the Laboratoire de Mod\'elisation et Contrôle de l'Universit\'e Joseph Fourier for numerous advices and helpful discussions.  
\appendix 
 
\section{Appendix : Plant modeling and simulation process} 
 
We first explain how the evolution of random variables is 
represented over the year. We then describe three important models 
of power generation units (thermal, hydro and demand side 
management contracts called EJP). We then briefly comment the 
resolution of the different dual subproblems. In subsection 
\ref{simulationprocess} the way of determining a production 
schedule is outlined.

\subsection{Prediction of random events and global problem} 
 
Several stochastic optimization problems could be envisaged to 
model the problem of yearly electric generation management (e.g. 
\cite{severalsto1,severalsto2,severalsto3}). The possible 
evolutions of the realization of the random variables are 
represented by a markov chain  (a set of scenarios organized in a 
tree). Each node of the tree corresponds to a one day period of 
time. A day is divided in $L$ hourly posts ($\mathcal{P}_n$ is the 
set of hourly posts of node $n$). At each node $n$ of the tree are 
attached a realization of the random variables we take into 
account : 
\begin{itemize} 
\item $\lambda_{n,p}$ is the price of electricity for post $p$. 
\item $\mathcal{D}_n(p)$ is the electricity demand for post $p$. 
\end{itemize} 
We will also use the following notations to describe the scenario tree: 
\begin{itemize} 
\item $\tau(n)$ is the time step associated to node $n$. 
\item $F(n)$ is the father node of node $n$. 
\item $S(n)$ is the set of son nodes of node $n$. 
\item $\pi_n$ is the probability to be at node $n$ ( $\displaystyle{\sum_{\tau(n)=t}} \; \pi_n=1 \;\; \forall \; t$) and $\pi_T(n)$ is the probability to go from $F(n)$ to $n$. 
\item $\mathcal{O}$ is the set of nodes ( $N= \card \; \mathcal{O}$). 
\item $T$ is the last time step and $\mathcal{O}_T$ is the set of leaves. 
\item $d_{n,p}$ is the duration of post $p$ at node $n$. 
\end{itemize} 
Other variables are attached to a node. They will be introduced in the description of the power generation units models. We give below the example of a scenario tree. In this example we have $F(1)=0$, $S(2)=\{4,5\},\ldots$ 
\begin{figure}[H] 
\begin{center} 
\includegraphics[angle=0, width=10cm,height=5cm]{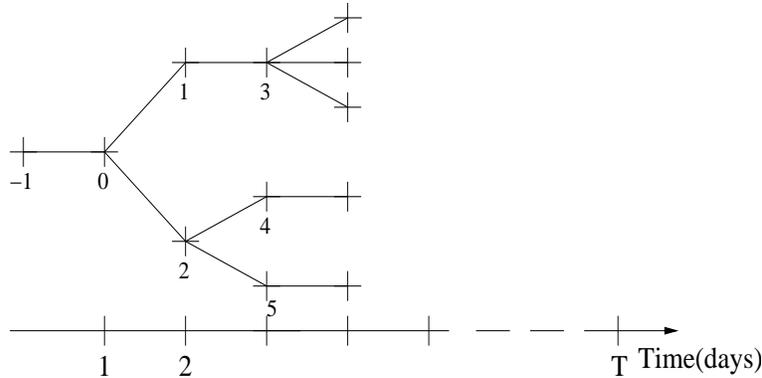} 
\end{center} 
\caption{Tree representing the different scenarios.} 
\end{figure} 
A scenario in the tree is thus a path from the root node $0$ to a leaf node. The construction of the tree is based on aggregation 
procedures using historical data 
(demand, inflows,...). Those problems are discussed in \cite{carpentier1},\cite{carpentier2}.\\ 
Given this representation of random events, the global problem of 
yearly power management scheduling consists in minimizing the 
average generation cost over the random tree while satisfying the 
demand constraint and the operating constraints of the generation 
units. It expresses as \eqref{base}. The reader should be aware 
that the solutions of this problem are indexed by the nodes of the 
scenario tree. If the scenario that occurs is represented in the 
tree we will have a generation schedule for this scenario. For a 
generation schedule that is not represented in the tree see 
subsection \ref{simulationprocess}. 
 
\subsection{Modelling of power generation units} \label{unitmodel} 
 
Three kinds of generation units are modelled : the thermal, hydro and EJP units. 
 
\paragraph{Thermal Units.} Let $\ell$ be the index of a thermal unit. The thermal units are completely described by two characteristics: 
\begin{enumerate} 
\item The generation levels :$P_{n,p}^\ell$ must remain between 
$P_{n,p}^{\ell,\min}=0$ and $P_{n,p}^{\ell,\max}=\tau_{f,n}^\ell 
\, \tau_{T,n}^\ell \, P_{\max}^\ell \, d_{n,p}$ with 
\begin{itemize} 
\item $\tau_{f,n}^\ell$ the random unavailability rate, realization of a random variable $\tau_\ell$. 
\item $\tau_{T,n}^\ell$ the programmed unavailability rate 
(deterministic). 
\item $P_{\max}^\ell$ the maximal power (in MW) of unit $\ell$. 
\end{itemize} 
\item The generation costs : $\forall n \in \mathcal{O}$ the thermal costs are 
$$ 
C_{n,p}^{\ell}(x_{n,p}^\ell,u_{n,p}^\ell)=c_\ell \, u_{n,p}^\ell, 
$$ 
where $c_\ell$ is the unitary production cost for unit $\ell$. 
\end{enumerate} 
 
\paragraph{Hydro Units.} An hydro valley is a set of interconnected plants and reservoirs with natural inflows in each reservoir. The state variables are the contents of the reservoir and the command variables are the discharge of the turbines and the water poured out of the reservoirs. In this yearly model, the total hydro production is aggregated in two different reservoirs non-interconnected. They both represent the total production capacity of the hydro generation units. The constraints are of two kinds: bounds on the volume of each reservoir, on the discharged water, and flow balance equations at each reservoir. Let $\ell$ be 
an hydro plant. We will use the following notations: 
\begin{itemize} 
\item $x_{n}^\ell$ is the content (in MWh) of the unique reservoir associated to $\ell$ at the beginning of time step $\tau(n)$. $x_{\min}^\ell$ and $x_{\max}^\ell$ are the lower and upper bounds (in MWh) on reservoir $\ell$ level. 
\item $\tau^\ell_{H,n}$ is the programmed unavailability rate (deterministic) 
\item $a_{n,p}^\ell$ is the natural inflow (in MWh) in reservoir $\ell$ for 
node $n$, post $p$. $a_n^\ell=\sum_{p \in \mathcal{P}_n} \; 
a_{n,p}^\ell$. 
\item $P_{\max}^\ell$ is the maximal power of hydro plant $\ell$. 
\item $dev_{n,p}^\ell$ is  the amount of water (in MWh) poured outside of reservoir $\ell$ and 
$v_{n,p}^\ell$ is the discharge of plant $\ell$ (in MWh). 
$v_{n,p}^{\ell,\min}=0$ and  $v_{n,p}^{\ell,\max}$ are the lower 
and upper bounds on the discharged power. Thus the command 
variables for the hydro subproblems are 
${(u_{n,p}^\ell)}_{n,p}={(v_{n,p}^\ell,dev_{n,p}^\ell)}_{n,p}$. 
\item $V_H^\ell(x)$ is the value of the water stock $x$ of reservoir $\ell$ at the last time step $T$. The operating costs for all reservoirs are null for $\tau(n)<T$. Hence, only the value of the water stock at time $T$ is taken into account. 
\end{itemize} 
If ${\mathcal{L}}_H$ is the set of hydro units, the hydro problem thus consists in minimizing 
$$ 
-\sum_{\ell \in \mathcal{L}_H} \; \sum_{n \in O_T}  \pi_n 
\;V_H^\ell(x_n^\ell+a_n^\ell-\sum_{p \in \mathcal{P}_n} \; 
v^\ell_{n,p}+ dev^\ell_{n,p}), 
$$ 
under the constraints 
\begin{equation} \label{hydrauconstraint} 
\left\{ 
\begin{array}{ll} 
x_n^\ell=x_{F(n)}^\ell +  \displaystyle{\sum_{p \in \mathcal{P}_n}} \; (a^\ell_{F(n),p}-v^\ell_{F(n),p}-dev^\ell_{F(n),p}) &\forall \;\; (n,\ell) \\ 
x_{min}^\ell \leq x_n^\ell +  \displaystyle{\sum_{p \in \mathcal{P}_n}} \; 
(a^\ell_{n,p}-v^\ell_{n,p}  -dev^\ell_{n,p}) \leq x_{\max}^\ell &\forall \;\; n \in \mathcal{O_T}\\ 
0 \leq v^\ell_{n,p} \leq  \displaystyle  {\tau}^\ell_{H,n} \, P_{max}^\ell \, d_{n,p} &  \forall \;\; (n,p,\ell)\\ 
0 \leq dev^\ell_{n,p} & \forall \;\; (n,p,\ell) \\ 
\displaystyle x_{min}^\ell \leq x_n^\ell \leq x_{max}^\ell & \forall \;\; (n,\ell). 
\end{array} 
\right. 
\end{equation} 
 
\paragraph{EJP contracts.} 
An EJP contract $\ell$ is represented by a production unit with the following features: 
\begin{itemize} 
\item $J_{\ell}$ is the total number of days the contract can be used. 
\item Each day, either the contract is used all day long or it is not 
used. A command variable $t_n^{\ell}$ defined for all node $n$ permits 
to know whether the contract $\ell$ is used at node $n$ 
($t_n^{\ell}=1$) or not ($t_n^{\ell}=0$). 
\item $s_n^{\ell}$ is the stock (in days) still available on the 
contract $\ell$ for node $n$ at the beginning of time step $\tau(n)$. 
\item $V_J^\ell(.)$ is the function defining the value of an EJP stock 
at the last time step $T$ for contract $\ell$. 
\item The power linked to contract $\ell$ is $P_J^\ell$. 
\end{itemize} 
Given a starting stock $J_{\ell}$ on EJP contract $\ell$ we maximize the value of the EJP stock at time step $T$ which yields to the following problem: 
\begin{equation} \label{EJPpb} 
\left\{ 
\begin{array}{ll} 
\min \;\; - \displaystyle{\sum_{\ell \in \mathcal{L}_J} \sum_{n \in \mathcal{O}_T}}  \; \pi_n V_J^\ell (s_n^\ell - t_n^\ell)&\\ 
s_n^\ell=J^\ell & \mbox{for the root node.}\\ 
t_n^\ell(t_n^\ell-1)=0 & \forall \; n  \in \mathcal{O}\\ 
s_{F(n)}^\ell -t_{F(n)}^\ell=s_n^\ell & \forall \; (n,\ell)\\ 
s_n^\ell \geq 0 & \forall \; (n,\ell)\\ 
s_{n}^\ell -t_{n}^\ell \geq 0 & \forall \; n  \in \mathcal{O}_T, 
\; \forall \ell, 
\end{array} 
\right. 
\end{equation} 
if $\mathcal{L}_J$ is the set of EJP contracts. Thus if $\ell \in 
\mathcal{L}_J$, the EJP command $u_{n,p}^\ell$ is given by $t_n \, 
P_J^\ell\,d_{n,p}$. To be complete, we should take into account 
nuclear power plants. However, this would not change much things 
as they can be modelled in a similar way (see \cite{nuclear4}). 
Using this modelling for the generation units, the constraint of 
satisfaction of the demand writes 
$\mathcal{D}_n(p)=\displaystyle{\sum_{\ell \in \mathcal{L}}} \; 
P_{n,p}^{\ell}(x_{n,p}^\ell,u_{n,p}^\ell)=\displaystyle{\sum_{\ell 
\in \mathcal{L}}} \; u_{n,p}^\ell$. 
 
\subsection{Solving the different dual subproblems} 
 
We briefly detail the first step of the space decomposition 
algorithm given in subsection \ref{spacedec}. The nominal thermal 
dual subproblem has an evident solution and the hydro subproblem 
is a linear optimization problem of big size (around 36 000 
variables in our case) solved using interior point methods. The 
EJP problem is an NP complete, non convex optimization problem 
solved using dynamic stochastic programming. We know, for every 
contract $\ell$, the values at the last time step $T$ of all 
possible values of the EJP stock $x$ : $V_J^\ell(x,T)$. We deduce, 
using HJB equations, (backward phase) the bellman values 
$V^\ell(x,n)$ for all contract $\ell$ and all node $n$: 
\begin{equation} \label{EJPsolve} 
V^\ell(x,n)=\left\{ 
\begin{array}{l} 
\max \; (\sum_{m \in S(n)} \; \pi_T(m) ( V^\ell(x - t_n^\ell,m) + 
\sum_{p \in \mathcal{P}_n} \; 
d_{n,p} \lambda_{n,p} P_J^\ell t_n^\ell ))\\ 
x-t_n^\ell \geq 0, \;\;\; t_n^\ell \in \{0,1\}. 
\end{array} 
\right.\\ 
\end{equation} 
Knowing the stock of every contract $\ell$ at the beginning of the 
year, the forward phase consists in deducing the optimal EJP 
commands to apply using the bellman values computed in the 
backward phase. Both kinds of methods used (dynamic stochastic 
programming and interior points methods) have a complexity that 
depends on the dimension of the dual space . A known drawback of 
dynamic programming is that its complexity grows exponentially 
with the state variable dimension. 
 
\subsection{The simulation process} \label{simulationprocess} 
 
The resolution of ($P$) provides optimal marginal prices 
$\lambda_{n,p}^*$ that are useful to elaborate a strategy which, 
for any realization of the random variables on the time period, 
will allow the computation of a generation schedule. This strategy 
has the form of Bellman functions and are computed with the 
following version of the Bellman principle. Let a reserve $\ell$ 
be given. For the last time step $T$, the Bellman function 
$V^\ell(x,n)$ is known for each stock level $x$ and node $n$ of 
the leaves. Between two grid points, the value is supposed to be 
linear. The Bellman values $V^\ell(x,n)$ for all the nodes $n$ of 
the tree and all stock step $x$ are computed using the following 
recursive formula: 
\begin{equation} \label{bellhyd} 
\left\{ 
\begin{array}{l} 
V^\ell(x,n)    = \displaystyle \max ( \sum_{m \in S(n)} \pi_{T}(m) 
( V^\ell(x + a_n^\ell-\sum_{p \in \mathcal{P}_n} \, v_{n,p}^\ell+dev_{n,p}^\ell,m)\\ 
\hspace*{3cm} + \displaystyle \sum_{p \in \mathcal{P}_n} 
\lambda_{n,p}^* \; v_{n,p}^\ell)) \\ 
x + a_n^\ell -x_{\max}^\ell \leq \displaystyle \sum_{p \in \mathcal{P}_n} \, v_{n,p}^\ell + dev_{n,p}^\ell \leq x + a_n^\ell\\ 
0 \leq v_{n,p}^\ell \leq v_{n,p}^{\ell,\max}, \quad 0 \leq dev_{n,p}^\ell 
\end{array} 
\right. 
\end{equation} 
if $\ell$ stands for an hydro reservoir and 
\begin{equation} \label{bellEJP} 
\left\{ 
\begin{array}{l} 
V^\ell(x,n)    = \max ( \displaystyle{\sum_{m \in S(n)}} 
\pi_{T}(m) ( V^\ell(x - t_n^\ell,m) + \displaystyle {\sum_{p \in 
\mathcal{P}_n}} 
{d_{n,p}} P_{J}^\ell \lambda_{n,p}^* t_n^\ell )) \\ 
t_n^\ell(t_n^\ell -1)=0, \quad x - t_n^\ell \geq 0 
\end{array} 
\right. 
\end{equation} 
if $\ell$ stands for an EJP contract. The Bellman function for time step $t$ is then given by: 
$ 
V^\ell(x,t)= \sum_{n \in \; t} \; \pi(n) V^\ell(x,n)$ for stock $x$, unit $\ell$. 
 
This algorithm is only a stochastic dynamic programming (SDP) performed on marginal values. Once those 
Bellman functions are computed, it is possible to perform a Monte-Carlo simulation of the generation 
scheduling using $\delta_x V^\ell(x,t)$ as a "fuel cost" of the energy kept in the reserve $\ell$. 
 
\bibliographystyle{unsrt} 
\bibliography{ValueAtRisk} 
\end{document}